# DYNAMIC IMPORTANCE SAMPLING FOR UNIFORMLY RECURRENT MARKOV CHAINS

By Paul Dupuis[1] and Hui Wang[2]

*Brown University*

Importance sampling is a variance reduction technique for efficient estimation of rare-event probabilities by Monte Carlo. In standard importance sampling schemes, the system is simulated using an a priori fixed change of measure suggested by a large deviation lower bound analysis. Recent work, however, has suggested that such schemes do not work well in many situations. In this paper we consider dynamic importance sampling in the setting of uniformly recurrent Markov chains. By "dynamic" we mean that in the course of a single simulation, the change of measure can depend on the outcome of the simulation up till that time. Based on a control-theoretic approach to large deviations, the existence of asymptotically optimal dynamic schemes is demonstrated in great generality. The implementation of the dynamic schemes is carried out with the help of a limiting Bellman equation. Numerical examples are presented to contrast the dynamic and standard schemes.

**1. Introduction.** Among variance reduction techniques for efficient Monte Carlo simulation is importance sampling, in which the data is generated using a probability distribution different from the true underlying distribution. It can be especially effective when applied to the estimation of expectations that are largely determined by rare events. To demonstrate the difficulty involved in simulating rare events by naive Monte Carlo, we consider a simple example. Let $X$ be a random variable taking values in $\mathbb{R}^d$, and suppose we are interested in estimating $p = P\{X \in A\}$ for some Borel set $A \subset \mathbb{R}^d$. To this end, a sequence of independent and identically distributed

Received June 2003; revised November 2003.

[1]Supported in part by the NSF Grants DMS-00-72004, ECS-99-79250 and Army Research Office Grants DAAD19-00-1-0549 and DAAD19-02-1-0425.

[2]Supported in part by NSF Grant DMS-01-03669.

(i.i.d.) copies $X_0, X_1, \ldots$ of $X$ are generated. With $I_k \doteq \mathbb{1}_{\{X_k \in A\}}$, an unbiased estimate for $p$ based on the first $K$ samples is just the sample mean: $Q_K \doteq (I_0 + I_1 + \cdots + I_{K-1})/K$. The relative error associated with this estimator is

$$\text{relative error} \doteq \frac{\text{standard deviation of } Q_K}{\text{mean of } Q_K} = \frac{\sqrt{p - p^2}}{p} \cdot \frac{1}{\sqrt{K}}.$$

Since $\sqrt{p - p^2}/p \to \infty$ as $p$ tends 0, a large sample size $K$ is required for the estimator $Q_K$ to achieve a reasonable relative error bound. For example, if $p = 10^{-8}$, ten billion samples are required to achieve a relative error bound of 10%.

The basic idea of importance sampling is as follows. Suppose that $X$ has distribution $\theta$, and consider an alternative sampling distribution $\tau$. It is required that $\theta$ be absolutely continuous with respect to $\tau$, so that the Radon–Nikodym derivative $f(x) \doteq (d\theta/d\tau)(x)$ exists. Independent and identically distributed samples $\bar{X}_0, \bar{X}_1, \ldots$ with distribution $\tau$ are generated. Form the estimate

$$\bar{Q}_K \doteq \frac{1}{K} \sum_{k=0}^{K-1} f(\bar{X}_k) \mathbb{1}_{\{\bar{X}_k \in A\}}$$

in lieu of $Q_K$. It is easy to check that $\bar{Q}_K$ is an unbiased estimate of $p$, with a rate of convergence determined by

$$\text{var}\,[f(\bar{X}_0)\mathbb{1}_{\{\bar{X}_0 \in A\}}] = \int_{\mathbb{R}} \mathbb{1}_{\{x \in A\}} f(x) \theta(dx) - p^2.$$

The optimization of this quantity over all possible $\tau$ is inappropriate. Indeed, taking $f(x) = p^{-1}\mathbb{1}_{\{x \in A\}}$ (i.e., $\tau$ is the conditional distribution of $X$ given $X \in A$), the variance becomes 0, but this change of measure requires the knowledge of the unknown parameter $p$. Instead, one typically seeks to minimize over parameterized families of alternative sampling distributions.

When the distribution of $X$ is connected to a large deviations problem, a standard heuristic is that the change of measure used to prove the large deviation lower bound should be a good (perhaps nearly optimal) distribution to use for the purposes of importance sampling. The first result of this type was given by Siegmund [34]. The basic idea was subsequently investigated in many contexts, and a small selection of the literally hundreds of papers on the topics is [[1, 2, 3, 7, 8, 9], [11, 12, 15, 17, 18, 20, 24, 25, 29, 30, 33]]. Necessary and sufficient conditions under which a prescribed scheme is asymptotically optimal are discussed in [10, 31, 32], while [21] gives a survey of rare-event simulation.

The validity of the heuristic, however, was challenged in [19]. Counterexamples were constructed to show that, under some very common settings, the change of measure suggested by large deviations leads to importance sampling schemes with very poor properties.

In order to explain these counterexamples, and more importantly, to find asymptotically optimal importance sampling algorithms in great generality, [16] introduces a dynamic importance sampling scheme and shows its asymptotic optimality in the setup of i.i.d. random variables (Cramér's theorem). The key observation is that *many* changes of measure are suggested by the large deviation lower bound analysis, and one must consider this larger class if one hopes to identify importance sampling schemes that work well in general. This leads to the development of schemes where the sampling distribution is dynamic (or, "adaptive") in the sense that the change of measure in the course of a single simulation can depend on the outcome of the simulation up till that time. For this reason, we also call such schemes *adaptive importance sampling* schemes.

The present paper analyzes the estimation of rare-event probabilities associated with uniformly recurrent Markov chains. More precisely, let $\{Y_j, j \in \mathbb{N}_0\}$ be a uniformly recurrent Markov chain taking values in a Polish space $\mathcal{S}$, and let $g : \mathcal{S} \to \mathbb{R}^d$ be a bounded measurable function. Define $S_n \doteq g(Y_0) + g(Y_1) + \cdots + g(Y_{n-1})$. The probability of interest is $P\{S_n/n \in A\}$ for a Borel set $A \subset \mathbb{R}^d$ and $n$ large. An asymptotic optimality result for traditional importance sampling is available in the one-dimensional case, $d = 1$, under the assumption which implies that the set $A$ is within a half interval that does not contain the expectation of $g$ under the invariant distribution [9]. A "dissection" approach was introduced for the high-dimensional case [9]. This approach was later on applied to Markov additive sequences [11], and was also implicitly used in [19]. This dissection approach requires that one appropriately partition the set $A$ into a finite number of subsets, and that a (possibly different) change of measure be applied to efficiently estimate the probability of each individual subset. However, there is no constructive way to obtain a suitable partition in general.

In this paper we develop adaptive importance sampling schemes for uniformly recurrent Markov chains. The existence of asymptotically optimal adaptive schemes is demonstrated for arbitrary dimension $d$, under very mild conditions on the set $A$. It turns out that one must study the asymptotics of a small noise stochastic *game* in order to analyze the optimality of importance sampling schemes. The distinction between the change of measures used in traditional importance sampling and adaptive importance sampling amounts, in control terminology, to the difference between "open-loop" and "feedback" controls. However, open loop controls are usually not optimal in the setting of stochastic games, except for very special cases. For this reason, the traditional importance sampling will not be asymptotically optimal in general. Our analysis indicates that the adaptive scheme also works for estimating functionals (other than probabilities) largely determined by rare events.

The paper is organized as follows. The setting of the problem is introduced in Section 2, with a brief description of the large deviations principle for uniformly recurrent Markov chains. We also give the definition of asymptotic optimality in this section. In Section 3 we show that adaptive importance sampling schemes designed to minimize the second moment are asymptotically optimal. Section 4 discusses an alternative formal PDE approach to the adaptive scheme, and describes a method for the construction of an asymptotically optimal adaptive scheme that does not directly depend on the large deviation parameter $n$. Numerical examples are presented in Section 4.3. Certain technical proofs are deferred to the appendices to ease exposition.

## 2. Problem setup and background.

2.1. *Problem setup.* Let $Y = \{Y_j, j \in \mathbb{N}_0\}$ be a time-homogeneous Markov chain taking values in a Polish space $\mathcal{S}$, with transition probability kernel

$$p(x, dy) = P\{Y_{j+1} \in dy | Y_j = x\}.$$

Let $g: \mathcal{S} \to \mathbb{R}^d$ be a bounded Borel-measurable function, and define

$$S_n \doteq g(Y_0) + g(Y_1) + \cdots + g(Y_{n-1}).$$

For an arbitrary Borel set $A \subset \mathbb{R}^d$, we wish to estimate

$$p_n \doteq P\{S_n/n \in A\}.$$

Throughout the paper we will make use of the following *uniform recurrency* assumption.

CONDITION 2.1. There exists a probability measure $\nu_p$ on $\mathcal{S}$, an integer $m_0 \in \mathbb{N}$ and a pair of strictly positive real numbers $a, b$ such that

$$a\nu_p(B) \leq p^{(m_0)}(x, B) \leq b\nu_p(B)$$

for all $x \in \mathcal{S}$ and Borel sets $B$. Here $p^{(m)}$ denotes the $m$-step probability transition kernel.

For example, an irreducible Markov chain with a finite state space is always uniformly recurrent.

The large deviation principle for a uniformly recurrent Markov chain is well known. It asserts that $\{S_n/n\}$ satisfies the large deviation principle with a convex rate function $L: \mathbb{R}^d \to [0, \infty]$. The identification of $L$ is deferred to the next section. We will impose the following assumption throughout the paper.

CONDITION 2.2. The Borel set $A \subset \mathbb{R}^d$ satisfies the condition

$$\inf_{\beta \in \bar{A}} L(\beta) = \inf_{\beta \in A^\circ} L(\beta).$$

Under Conditions 2.1 and 2.2, we have the large deviations approximation

$$\lim_{n \to \infty} \frac{1}{n} \log P\{S_n/n \in A\} = - \inf_{\beta \in A} L(\beta).$$

REMARK 2.1. The uniform recurrency assumption (Condition 2.1) is convenient to work with. It includes the important case of irreducible finite state Markov chains, and generalizes the results in [16] where i.i.d. sequences were considered. However, this strong recurrency assumption also excludes many important Markov chains. One difficulty in extending the present results to more general Markov chains is that the uniform positivity and boundedness of the eigenfunctions (see Section 2.2) may not be preserved [26, 27]. It is clear that generalization in this direction will require a much more involved analysis.

2.2. *LDP for a uniformly recurrent Markov chain.* In this section we discuss two different approaches to the identification of the rate function $L$. The first approach suggests a parameterized family of change of measures (see Remark 2.2) that will be used later on to build importance sampling schemes. The second approach identifies the rate function $L$ in terms of relative entropy, and will be used in the analysis of the asymptotic optimality of adaptive schemes.

The first approach is based on a generalized Perron–Frobenius theorem. Fix any $\alpha \in \mathbb{R}^d$. Then by [22], the nonnegative kernel

$$\exp\{\langle \alpha, g(y) \rangle\} p(x, dy)$$

admits a unique real eigenvalue $\exp\{H(\alpha)\}$ and a unique (up to a multiplicative constant) eigenfunction $r(x; \alpha)$ in the sense that, for every $x \in \mathcal{S}$,

$$(2.1) \qquad \int_{\mathcal{S}} e^{\langle \alpha, g(y) \rangle} r(y; \alpha) p(x, dy) = e^{H(\alpha)} r(x; \alpha),$$

and with the following properties. $H(\alpha)$ is an analytic, strictly convex function of $\alpha \in \mathbb{R}^d$ with $H(0) = 0$, and there exist $0 < c_\alpha < C_\alpha < \infty$ such that

$$(2.2) \qquad c_\alpha \leq r(x; \alpha) \leq C_\alpha \qquad \forall\, x \in \mathcal{S}.$$

The paper [22] also shows that the rate function of the large deviation principle for $\{S_n/n\}$ is the convex conjugate of $H$, that is,

$$(2.3) \qquad L(\beta) = \sup_{\alpha \in \mathbb{R}^d} [\langle \alpha, \beta \rangle - H(\alpha)].$$

Note that in the special case when the Markov chain $Y$ is an i.i.d. sequence, $H(\alpha)$ is the logarithm moment generating function of $g(Y_j)$ and $r(x;\alpha) \equiv 1$. Therefore, this result generalizes the classical Cramér's theorem, at least for bounded i.i.d. random variables. For the case when $Y$ is an irreducible Markov chain with finite state space, $\exp\{H(\alpha)\}$ is just the maximal eigenvalue of the irreducible nonnegative matrix $\exp\{\langle \alpha, g(y)\rangle\}p(x,dy)$, and $r(\cdot;\alpha)$ is the associated right eigenvector.

REMARK 2.2. It is not difficult to see that, thanks to (2.1), for each $\alpha \in \mathbb{R}^d$,

$$\exp\{\langle \alpha, g(y)\rangle - H(\alpha)\} \cdot \frac{r(y;\alpha)}{r(x;\alpha)} \cdot p(x,dy)$$

defines a probability transition kernel.

Another approach is the weak convergence methodology which utilizes a stochastic control representation for certain exponential integrals [14]. It first identifies the large deviations rate function for the empirical measure of the Markov chain in the $\tau$-topology, then uses contraction principle to obtain the rate function for $\{S_n/n\}$. We will need the following definitions.

For an arbitrary Polish space $\mathcal{Z}$, we denote by $\mathcal{P}(\mathcal{Z})$ the collection of all probability measures on space $(\mathcal{Z}, \mathcal{B}(\mathcal{Z}))$. For a pair of probability measures $\gamma, \mu \in \mathcal{P}(\mathcal{Z})$, the *relative entropy* of $\gamma$ with respect to $\mu$ is defined as

$$R(\gamma\|\mu) \doteq \begin{cases} \int_{\mathcal{Z}} \log \frac{d\gamma}{d\mu} d\gamma, & \text{if } \gamma \ll \mu, \\ \infty, & \text{otherwise.} \end{cases}$$

Given a probability transition kernel $q(x,dy)$ on space $\mathcal{Z}$, we define $\mu q \in \mathcal{P}(\mathcal{Z})$, $\mu \otimes q \in \mathcal{P}(\mathcal{Z} \times \mathcal{Z})$ by

$$\mu q(B) \doteq \int_{\mathcal{Z}} q(x,B)\,\mu(dx),$$

$$(\mu \otimes q)(D \times B) \doteq \int_{D \times B} \mu(dx)q(x,dy) = \int_D q(x,B)\,\mu(dx)$$

for all Borel sets $D, B \subset \mathcal{Z}$. The collection of all probability transition kernels on $\mathcal{Z}$ is denoted by $\mathcal{T}(\mathcal{Z})$.

The weak convergence approach identifies the rate function for $\{S_n/n\}$ in terms of relative entropy:

$$L(\beta) = \inf \Big\{ R(\mu \otimes q \| \mu \otimes p) : \mu \in \mathcal{P}(\mathcal{S}),$$

(2.4)
$$q \in \mathcal{T}(\mathcal{S}), \mu q = \mu, \int_{\mathcal{S}} g\,d\mu = \beta \Big\}.$$

The validity of the representation (2.4) is implied by the results in [14], Chapters 8 and 9, where the large deviation principle of the empirical measures associated with Markov chains are studied under weaker assumptions.

For future reference, we summarize the preceding discussion into the following proposition. The only part that has not been mentioned is the superlinearity of the rate function $L$, which is an easy consequence of (2.3) and the finiteness of $H$ ([14], Lemma 6.2.3(c)).

PROPOSITION 2.1. *Under Condition* 2.1, *the sequence* $\{S_n/n\}$ *satisfies the large deviation principle with rate function* $L$, *which is given by* (2.3) *and* (2.4). *Moreover, the rate function* $L$ *is convex, lower-semicontinuous and superlinear in the sense that*

$$\lim_{N \to \infty} \inf_{\{\beta \in \mathbb{R}^d : \|\beta\| \geq N\}} \frac{L(\beta)}{\|\beta\|} = \infty.$$

*In particular,* $L$ *has compact level sets.*

2.3. *Asymptotic optimality.* In this section we define *asymptotic optimality* for an importance sampling scheme.

Consider a probability space $(\Omega, \mathcal{F}, P)$ and a family of events $\{A_n\}$ with

$$\lim_{n \to \infty} \frac{1}{n} \log P\{A_n\} = -\gamma,$$

for some $\gamma \geq 0$. A general formulation of importance sampling for this problem can be described as follows. In order to estimate $P\{A_n\}$, a generic random variable $\bar{Z}_n$ is constructed such that $P\{A_n\} = E\bar{Z}_n$. Independent replications $(\bar{Z}_n^0, \bar{Z}_n^1, \ldots, \bar{Z}_n^{K-1})$ of $\bar{Z}_n$ are then generated, and we obtain an estimator by averaging

$$\bar{Q}_n^K \doteq \frac{\bar{Z}_n^0 + \bar{Z}_n^1 + \cdots + \bar{Z}_n^{K-1}}{K}.$$

The estimator is unbiased, that is, $E\bar{Q}_n^K = P\{A_n\}$. The rate of convergence associated with this estimator is determined by the variance of the summands, or equivalently, their second moment $E[(\bar{Z}_n)^2]$. The smaller the second moment, the faster the convergence, whence the smaller sample size $K$ required. However, it follows from Jensen's inequality that

$$\limsup_{n \to \infty} -\frac{1}{n} \log E[(\bar{Z}_n)^2] \leq \lim_{n \to \infty} -\frac{1}{n} \log (E\bar{Z}_n)^2 = 2\gamma.$$

The estimator $\bar{Q}_n^K$ is said to be *asymptotically optimal* if

$$\lim_{n \to \infty} -\frac{1}{n} \log E[(\bar{Z}_n)^2] = 2\gamma.$$

REMARK 2.3. Since the performance of the estimator $\bar{Q}_n^K$ is completely determined by the second moment of its generic, i.i.d. building block $\bar{Z}_n^k$, we will drop the superscript $k$ hereafter. Note that $n$ does *not* stand for sample size, but for the large deviation parameter.

**3. Statement of the main result.** The adaptive importance sampling scheme we consider dynamically selects the change of measure (or the parameter $\alpha$) in the form suggested by Remark 2.2, according to the sample history. Naturally, the scheme is closely related to a control problem. Let the control $\alpha^n = \{\alpha_j^n(\cdot,\cdot), j = 1, \ldots, n-1\}$ be given, where each $\alpha_j^n : \mathcal{S} \times \mathbb{R}^d \to \mathbb{R}^d$ is a Borel-measurable function. Then the state dynamics are governed by

$$\bar{S}_j^n \doteq \sum_{i=0}^{j-1} g(\bar{Y}_i^n), \qquad j = 0, 1, \ldots, n.$$

Here we set $\bar{Y}_0 = Y_0 \equiv y_0$, and for $j \geq 1$, $\bar{Y}_j^n$ is conditionally distributed, given $\{\bar{Y}_i^n, i = 0, 1, \ldots, j-1\}$, according to

$$v_j^n(dy) = \exp\{\langle \alpha_j^n, g(y) \rangle - H(\alpha_j^n)\} \cdot \frac{r(y; \alpha_j^n)}{r(\bar{Y}_{j-1}^n; \alpha_j^n)} \cdot p(\bar{Y}_{j-1}^n, dy)$$

with (abusing notation a bit) $\alpha_j^n = \alpha_j^n(\bar{Y}_{j-1}^n, \bar{S}_j^n/n)$.

An unbiased estimator of $P\{S_n/n \in A\}$ is defined as the average of independent copies of

$$\bar{X}_n = \mathbb{1}_{\{\bar{S}_n^n/n \in A\}} \exp\left\{\sum_{j=1}^{n-1} (-\langle \alpha_j^n, g(\bar{Y}_j^n)\rangle + H(\alpha_j^n))\right\} \cdot \prod_{j=1}^{n-1} \frac{r(\bar{Y}_{j-1}^n; \alpha_j^n)}{r(\bar{Y}_j^n; \alpha_j^n)}.$$

Our goal is to minimize the second moment, hence the variance, of the summands $\bar{X}_n$ by judiciously choosing the control $\alpha^n$. Thus, we consider the value function defined by

$$V^n(y_0) \doteq \inf_{\alpha^n} E[\bar{X}_n^2]$$

$$= \inf_{\alpha^n} E\left[\mathbb{1}_{\{\bar{S}_n^n/n \in A\}} \exp\left\{\sum_{j=1}^{n-1}(-2\langle \alpha_j^n, g(\bar{Y}_j^n)\rangle + 2H(\alpha_j^n))\right\}\right.$$
$$\left. \times \prod_{j=1}^{n-1} \frac{r^2(\bar{Y}_{j-1}^n; \alpha_j^n)}{r^2(\bar{Y}_j^n; \alpha_j^n)}\right].$$

For convenience we write $V^n(y_0)$ as $V^n$ when no confusion is incurred. We also consider the log transform

$$W^n = -\frac{1}{n} \log V^n.$$

We have the following result, which asserts the existence of asymptotically optimal adaptive importance sampling schemes.

THEOREM 3.1. *Under Conditions 2.1 and 2.2, we have*
$$\lim_{n\to\infty} W^n = 2 \inf_{\beta \in A} L(\beta).$$

The detailed proof is deferred to Appendix A. It is worth pointing out that the construction of asymptotically optimal or nearly optimal adaptive schemes (i.e., selection of the control $\alpha^n$) is implied by a *dynamic programming equation* (DPE) appearing in the proof. Since the proof is rather lengthy and technical, it makes sense to give an outline and some intuitive discussion below, so that readers can proceed to the construction of the adaptive schemes (Section 4), without having to delve into the technical details of the proof.

*Outline and intuition of the proof.* Thanks to the discussion in Section 2.3, it suffices to show the lower bound

(3.1) $$\liminf_{n} W^n \geq 2 \inf_{\beta \in A} L(\beta).$$

The proof will utilize the DPE that is satisfied by $W^n$. In order to do so, we first extend the dynamics. Abusing notation a bit, for $x \in \mathbb{R}^d$, $y \in \mathcal{S}$ and $i \in \{0, 1, \ldots, n\}$, define the dynamics

$$\bar{S}_{i,j}^n = nx + \sum_{\ell=i}^{j-1} \bar{Y}_{i,\ell}^n, \qquad j = i, \ldots, n.$$

Here we set $\bar{Y}_{i,i} \equiv y$, and for $j \geq i+1$, $\bar{Y}_{i,j}^n$ is conditionally distributed, given $\{\bar{Y}_{i,\ell}, \ell = i, \ldots, j-1\}$, according to

$$v_{i,j}^n(dz) = \exp\{\langle \alpha_j^n, g(z) \rangle - H(\alpha_j^n)\} \cdot \frac{r(z; \alpha_j^n)}{r(\bar{Y}_{i,j-1}^n; \alpha_j^n)} p(\bar{Y}_{i,j-1}^n, dz),$$

where $\alpha_j^n = \alpha_j^n(\bar{Y}_{i,j-1}^n, \bar{S}_{i,j}^n/n)$. The original control problem corresponds to $x = 0$, $i = 0, y = y_0$. Define analogously

$$V^n(x, y; i)$$
$$\doteq \inf_{\alpha^n} E\left[ \mathbb{1}_{\{\bar{S}_{i,n}^n/n \in A\}} \exp\left\{ \sum_{j=i+1}^{n-1} (-2\langle \alpha_j^n, g(\bar{Y}_{i,j}^n) \rangle + 2H(\alpha_j^n)) \right\} \right.$$
$$\left. \times \prod_{j=i+1}^{n-1} \frac{r^2(\bar{Y}_{i,j-1}^n; \alpha_j^n)}{r^2(\bar{Y}_{i,j}^n; \alpha_j^n)} \right]$$

and its log transform

$$W^n(x, y; i) = -\frac{1}{n} \log V^n(x, y; i).$$

The terminal conditions are

$$V^n(x,y;n) = \mathbb{1}_A(x), \qquad W^n(x,y;n) = \infty \cdot \mathbb{1}_{A^c}(x).$$

Since it is inconvenient to study a problem with an $\infty$ terminal condition, we instead work with a mollified version of the control problem. Let $F : \mathbb{R}^d \to \mathbb{R}$ be an arbitrary bounded and Lipschitz continuous function. Suppose that $V_F^n$ is defined as $V^n$, save that the indicator function $\mathbb{1}_{\{\bar{S}_{i,n}^n/n \in A\}}$ is replaced by $\exp\{-2nF(\bar{S}_{i,n}^n/n)\}$. Similarly define

$$(3.2) \qquad W_F^n(x,y;i) \doteq -\frac{1}{n} \log V_F^n(x,y;i).$$

Since $V_F^n$ is the value function of a control problem, one can write down the DPE for $V_F^n$. Substituting (3.2) in this DPE, one obtains an equation for $W_F^n$; see (A.1). The proof of the desired inequality (3.1) is based on the analysis of this recursive equation for $W_F^n$.

The relative entropy representation for exponential integrals ([14], Proposition 1.4.2) states that

$$(3.3) \qquad -\log \int_{\mathcal{S}} e^{-f(x)} \mu(dx) = \inf_{\gamma \in \mathcal{P}(\mathcal{S})} \left[ R(\gamma \| \mu) + \int f \, d\gamma \right]$$

for all bounded and Borel measurable functions $f$. Applying this representation formula to the equation (A.1) for $W_F^n$, one obtains

$$(3.4) \quad \begin{aligned} W_F^n&(x,y;i) \\ &= \sup_{\alpha \in \mathbb{R}^d} \inf_{\gamma \in \mathcal{P}(\mathcal{S})} \bigg[ \int W_F^n\bigg(x + \frac{1}{n}g(y), z; i+1\bigg) \gamma(dz) \\ &\qquad + \frac{1}{n}\bigg( R(\gamma(\cdot)\|p(y,\cdot)) + \int \langle \alpha, g(z) \rangle \gamma(dz) - H(\alpha) \bigg) \\ &\qquad + \frac{1}{n} \int \log \frac{r(z;\alpha)}{r(y;\alpha)} \gamma(dz) \bigg]. \end{aligned}$$

This equation suggests that $W_F^n$ is the lower value of a discrete-time stochastic game. One of the two players of the game (the $\alpha$-player) selects the parameter $\alpha$, and is the weaker player. The other player (the $\gamma$-player) is the stronger player, and selects the distribution $\gamma$ that determines the evolution of the state. The right-hand side of (3.4) would take a simpler form if we could permute the sup and inf. However, this is not (in general) possible, since the last term

$$(3.5) \qquad \frac{1}{n} \int \log \frac{r(z;\alpha)}{r(y;\alpha)} \gamma(dz)$$

may not be concave in $\alpha$.

This difficulty is also the main distinction from the setting of Cramér's theorem where the Markov chain $Y$ reduces to an i.i.d. sequence of random variables. The latter case gives $r(x;\alpha) \equiv 1$ and the unpleasant term (3.5) disappears, whence the min/max theorem can be applied to convert the DPE of $W_F^n$ into a DPE associated with a *control* problem, which is much simpler to analyze than a game [16]. However, the interchange of sup and inf is not possible with (3.4) as written.

The key idea to overcome this difficulty and to obtain a lower bound for $W_F^n$ is as follows. Fix an integer $m$, and consider a variant of the game where the $\alpha$-player is constrained to policies such that $\alpha$ must be constant over time intervals of length $1/m$. This new game is even more favorable to the $\gamma$-player, whence it will have a smaller lower-value. Letting $n$ go to infinity, the lower value of the new game converges to a function $U_F^m$, and we expect

$$\liminf_{n \to \infty} W_F^n(x, y; \lfloor nk/m \rfloor) \geq U_F^m(x; k), \qquad k = 0, 1, \ldots, m.$$

A bonus of taking the limit is that the troubling terms (3.5), which can be interpreted as part of the running cost, cancel off, and it is not difficult to guess that $U_F^m$ should satisfy

$$(3.6) \quad U_F^m(x; k) = \sup_{\alpha \in \mathbb{R}^d} \inf_{\beta \in \mathbb{R}^d} \left[ U_F^m\left(x + \frac{1}{m}\beta; k+1\right) + \frac{1}{m}(L(\beta) + \langle \alpha, \beta \rangle - H(\alpha)) \right],$$

with terminal condition

$$(3.7) \qquad U_F^m(x; m) \doteq 2F(x).$$

In the proof, $U_F^m$ is in fact defined recursively through equations (3.6) and (3.7).

Equation (3.6) is much easier to analyze. Analogously to [16], one can show by a weak convergence argument that

$$(3.8) \qquad \liminf_{m \to \infty} U_F^m(x, 0) \geq 2 \inf_{\beta \in \mathbb{R}^d} \{L(\beta) + F(x + \beta)\},$$

which in turn implies

$$\liminf_{n \to \infty} W_F^n(x, y; 0) \geq 2 \inf_{\beta \in \mathbb{R}^d} \{L(\beta) + F(x + \beta)\}.$$

Letting $x = 0$ and the mollifier $F$ tend to $\infty \cdot \mathbb{1}_{A^c}$, one arrives at the desired inequality (3.1). □

The following result is useful in the identification of an optimal adaptive importance sampling scheme in Section 4.

COROLLARY 3.2. *Fix an arbitrary $x \in \mathbb{R}^d$, and a bounded Lipschitz continuous function $F : \mathbb{R}^d \to \mathbb{R}$. Assume Condition* 2.1, *and define $U_F^m$ recursively by* (3.6) *with the terminal condition* (3.7). *Then*

$$\lim_{m \to \infty} U_F^m(x; \lfloor tm \rfloor) = 2U_F(x, t) \qquad \forall t \in [0, 1],$$

*where*

(3.9) $$U_F(x, t) \doteq \inf_{\beta \in \mathbb{R}^d} \{(1-t)L(\beta) + F(x + (1-t)\beta)\}.$$

PROOF. We will show the equality for $t = 0$. The case with general $t \in [0, 1]$ is similar and thus omitted.

Thanks to (3.8), it suffices to prove

$$\limsup_{m \to \infty} U_F^m(x; 0) \le 2U_F(x, 0) = 2 \inf_{\beta \in \mathbb{R}^d} \{L(\beta) + F(x + \beta)\}.$$

Fix an arbitrary $\beta \in \mathbb{R}^d$. The recursive definition of $U_F^m$ (3.6) and (2.3) yield

$$U_F^m(x; k) \le \sup_{\alpha \in \mathbb{R}^d} \left[ U_F^m\left(x + \frac{1}{m}\beta; k+1\right) + \frac{1}{m}(L(\beta) + \langle \alpha, \beta \rangle - H(\alpha)) \right]$$
$$= U_F^m\left(x + \frac{1}{m}\beta; k+1\right) + \frac{2}{m}L(\beta).$$

Repeatedly applying this inequality for $k = 0, 1, \ldots, m-1$, we arrive at

$$U_F^m(x; 0) \le U_F^m(x + \beta; m) + 2L(\beta) = 2F(x + \beta) + 2L(\beta),$$

thanks to (3.7). This completes the proof. $\square$

## 4. Implementation issues and examples.

4.1. *The limit control problem and implementation issues.* Theorem 3.1 establishes the existence of asymptotically optimal adaptive sampling schemes. However, it does not explicitly discuss the construction of such schemes. On the other hand, the proof of the theorem implies that one approach of construction would be to solve, numerically if need be, the DPE (3.4) associated with $W_F^n$ ($W^n$ equals $W_F^n$ when $F = \infty \cdot \mathbb{1}_{A^c}$). However, this approach may not only require a lot of computation effort, but the resulting adaptive sampling control (i.e., control $\alpha^n$) will directly depend on $n$. In general, one would prefer schemes without this dependence.

An alternative approach is to consider the DPE associated with the limit problem of $U_F^m$ as $m$ tends to infinity. To this end, we rewrite (3.6) as

$$0 = \sup_{\alpha \in \mathbb{R}^d} \inf_{\beta \in \mathbb{R}^d} \left[ \triangle U_F^m + \frac{1}{m}(L(\beta) + \langle \alpha, \beta \rangle - H(\alpha)) \right],$$

where

$$\triangle U_F^m \doteq U_F^m\left(x + \frac{1}{m}\beta; k+1\right) - U_F^m(x; k).$$

Suppose that a $t$ subscript denotes the partial derivative with respect to $t$, and that an $x$ subscript denotes the vector of partials with respect to $x_i$, $i = 1, \ldots, d$. Since Corollary 3.2 (for $F$ bounded and Lipschitz continuous) asserts that

$$\lim_{m \to \infty} U_F^m(x; \lfloor tm \rfloor) = 2U_F(x; t),$$

we have formally the approximation

$$\triangle U_F^m \approx \left\langle \frac{1}{m}\beta, (2U_F)_x \right\rangle + \frac{1}{m}(2U_F)_t.$$

Substituting this back, we have

$$0 = \sup_{\alpha \in \mathbb{R}^d} \inf_{\beta \in \mathbb{R}^d} [\langle \beta, (2U_F)_x \rangle + (2U_F)_t + L(\beta) + \langle \alpha, \beta \rangle - H(\alpha)]$$
$$= (2U_F)_t + \sup_{\alpha \in \mathbb{R}^d} \inf_{\beta \in \mathbb{R}^d} [L(\beta) + \langle \alpha + (2U_F)_x, \beta \rangle - H(\alpha)].$$

Representing the infimum in terms of the Legendre transform $H$ of $L$ gives

$$0 = (2U_F)_t + \sup_{\alpha \in \mathbb{R}^d} [-H(-\alpha - (2U_F)_x) - H(\alpha)].$$

The strict convexity of $H$ implies that

(4.1) $$\alpha^*(x, t) = -(U_F)_x(x, t),$$

and that

(4.2) $$0 = (U_F)_t - H(-(U_F)_x).$$

Equation (4.1) identifies, at least formally, an optimal feedback control policy. However, this observation is not entirely satisfactory since $U_F$ does not usually have an explicit solution, and even if there is an exact formula for $U_F$, the partial derivatives may not be defined for all time and spatial points. In order to obtain a formal characterization of $\alpha^*$ that is more useful, we observe that, thanks to the definition (3.9) of $U_F$ and the convexity of $L$, $U_F$ is the value function of the deterministic control problem

$$U_F(x, t) = \inf_\phi \left[ \int_t^1 L(\dot\phi(s))\,ds + F(\phi(1)) \right],$$

where the infimum is over all absolutely continuous $\phi$ which satisfy $\phi(t) = x$. It is straightforward to see from this control problem that an optimal control at $(x, t)$ is the minimizer in (3.9), say $\beta^*(x, t)$, thanks to the convexity of $L$.

The standard dynamic programming argument implies that $U_F$ (in a weak sense) satisfies the DPE

$$0 = (U_F)_t + \inf_{a \in \mathbb{R}^d}[L(a) + \langle a, (U_F)_x \rangle] = (U_F)_t - H(-(U_F)_x),$$

which, not surprisingly, is just equation (4.2). The optimal control $\beta^*(x,t)$ is, at least formally, the minimizer in the DPE, or $\beta^*(x,t)$ and $-(U_F)_x(x,t)$ are conjugate. It follows that

$\alpha^*(x,t)$ is conjugate to the minimizer $\beta^*(x,t)$ in (3.9).

At points where $(U_F)_x(x,t)$ exists this definition gives $\alpha^*(x,t) = -(U_F)_x(x,t)$. At points where $(U_F)_x(x,t)$ does not exist there are multiple minimizing $\beta^*(x,t)$, and one should define $\alpha^*(x,t)$ through conjugacy in any Borel measurable way.

REMARK 4.1. The original (unmollified) problem corresponds to $F = \infty \cdot \mathbb{1}_{A^c}$. In this case,

(4.3) $\qquad \beta^*(x,t) \in \arg\min\{(1-t)L(\beta) : x + (1-t)\beta \in A\},$

and $\alpha^*(x,t)$ is its conjugate.

4.2. *Numerical examples.* We give two numerical examples in order to illustrate the asymptotic optimality of the adaptive schemes, in general, and the pitfalls of the traditional importance sampling schemes. The first example is concerned with a simple Markov chain with two states, while the second example studies a discrete time Markov chain embedded in a tandem Jackson network with finite buffers.

EXAMPLE 4.1. Consider a simple finite-state Markov chain $Y$ with state space $\mathcal{S} = \{1, -1\}$ and probability transition matrix

$$Q = \begin{bmatrix} p & 1-p \\ 1 & 0 \end{bmatrix} \doteq [Q(i,j)]_{2 \times 2},$$

for some constant $p \in (0,1)$. Define $g : \mathcal{S} \to \mathbb{R}$ by $g(x) \doteq x$, and $S_n = g(Y_0) + g(Y_1) + \cdots + g(Y_{n-1}) = Y_0 + Y_1 + \cdots + Y_{n-1}$.

Since $Y$ is an irreducible finite-state Markov chain, the eigenvalue $e^{H(\alpha)}$ and eigenfunction $r(\cdot; \alpha)$, as defined in (2.1) for $\alpha \in \mathbb{R}$, are just the maximal eigenvalue of the kernel $[e^{\alpha j}Q(i,j)]$ and the corresponding eigenvector, respectively. Simple algebra gives

$$H(\alpha) = \log \frac{pe^\alpha + \sqrt{p^2 e^{2\alpha} + 4(1-p)}}{2} \qquad \forall \alpha \in \mathbb{R},$$

which is a convex function with $H(0) = 0$, and an eigenvector

$$\begin{bmatrix} r(1;\alpha) \\ r(-1;\alpha) \end{bmatrix} = \begin{bmatrix} e^{H(\alpha)} \\ e^{\alpha} \end{bmatrix}.$$

Therefore, for any given $\alpha \in \mathbb{R}$, the corresponding change of measure is represented by the probability transition matrix

(4.4)
$$Q_\alpha = \left[ e^{\alpha j - H(\alpha)} \frac{r(j;\alpha)}{r(i;\alpha)} Q(i,j) \right]$$
$$= \begin{bmatrix} pe^{\alpha - H(\alpha)} & (1-p)e^{-2H(\alpha)} \\ 1 & 0 \end{bmatrix}.$$

Let $L$ be the convex conjugate of $H$. It is not difficult to check that $L(\beta) = \infty$ if $\beta < 0$ or $\beta > 1$, and that for $\beta \in (0,1)$,

$$L(\beta) = \sup_{\alpha \in \mathbb{R}} [\alpha\beta - H(\alpha)]$$
$$= \frac{\beta}{2} \log \frac{4(1-p)\beta^2}{p^2(1-\beta^2)} + \frac{1}{2} \log \frac{1-\beta}{1+\beta} - \frac{1}{2} \log(1-p)$$

with the minimizer

(4.5) $$\alpha^* \doteq \alpha^*(\beta) = \frac{1}{2} \log \frac{4(1-p)\beta^2}{p^2(1-\beta^2)}$$

and

$$L(0) = \lim_{\beta \downarrow 0} L(\beta) = -\tfrac{1}{2} \log(1-p),$$

$$L(1) = \lim_{\beta \uparrow 1} L(\beta) = -\log p.$$

Furthermore, $L(\beta) = 0$ if and only if $\beta = H'(0) = p/(2-p)$.

Assume $Y_0 \equiv 1$. We are interested in estimating $p_n \doteq P\{S_n/n \in A\}$ for the Borel set

$$A = (-\infty, a] \cup [b, \infty), \qquad 0 < a < H'(0) < b < 1.$$

In all the following discussion we take $p = 1/2$, $a = 1/6$, $b = 1/2$, which implies

(4.6) $$\inf_{\beta \in A} L(\beta) = L(b) < L(a).$$

We will compare the naive Monte Carlo simulation, traditional importance sampling and adaptive importance sampling schemes below.

The naive Monte Carlo simulation will simulate the Markov chain under the original transition probability kernel $Q$. One can also regard this as a special change of measure with the corresponding $\alpha = 0$. In this case, the estimate is just the sample mean of $K$ i.i.d. replications of $X_n = \mathbb{1}_{\{S_n/n \in A\}}$. Since the second moment of $X_n$ satisfies

$$\lim_{n \to \infty} -\frac{1}{n} \log E[(X_n)^2] = \lim_{n \to \infty} -\frac{1}{n} \log p_n$$
$$= \inf_{\beta \in A} L(\beta)$$
$$= L(b) < 2L(b),$$

the naive Monte Carlo sampling is not asymptotically optimal.

Thanks to (4.6), the traditional importance sampling will take $\beta^* = b$, and $\alpha^*$ is then defined by (4.5). The algorithm will generate a Markov chain $\tilde{Y}$ with probability transition matrix $Q_{\alpha^*}$ and $\tilde{Y}_0 \equiv 1$. Let

$$\tilde{S}_n \doteq \tilde{Y}_0 + \cdots + \tilde{Y}_{n-1}.$$

The estimate is the sample mean of $K$ i.i.d. replications of

$$\tilde{X}_n = \mathbb{1}_{\{\tilde{S}_n/n \in A\}} \prod_{j=1}^{n-1} e^{-\alpha^* \tilde{Y}_j + H(\alpha^*)} \cdot \prod_{j=1}^{n-1} \frac{r(\tilde{Y}_{j-1}; \alpha^*)}{r(\tilde{Y}_j; \alpha^*)}$$
$$= \mathbb{1}_{\{\tilde{S}_n/n \in A\}} e^{-\alpha^* \tilde{S}_n + nH(\alpha^*)} \cdot e^{\alpha^* \tilde{Y}_0 - H(\alpha^*)} \frac{r(\tilde{Y}_0; \alpha^*)}{r(\tilde{Y}_{n-1}; \alpha^*)}.$$

Since $r(\cdot; \alpha^*)$ is clearly bounded from above and bounded away from zero, it is not difficult to see that

$$\lim_{n \to \infty} -\frac{1}{n} \log E[(\tilde{X}_n)^2] = \lim_{n \to \infty} -\frac{1}{n} \log E[\mathbb{1}_{\{\tilde{S}_n/n \in A\}} e^{-2n(\alpha^* \tilde{S}_n/n - H(\alpha^*))}].$$

Simple computation yields that $\{\tilde{S}_n/n\}$ satisfies the large deviation principle with rate function $\tilde{L}(\beta) = L(\beta) + H(\alpha^*) - \alpha^* \beta$. Now one can apply the Varadhan's theorem ([14], Theorem 1.3.4) (with slight modification) to show

$$\lim_{n \to \infty} -\frac{1}{n} \log E[(\tilde{X}_n)^2] = \inf_{\beta \in A} [2\alpha^* \beta - 2H(\alpha^*) + \tilde{L}(\beta)]$$
$$= \inf_{\beta \in A} [\alpha^* \beta - H(\alpha^*) + L(\beta)].$$

In the configuration of this example, the infimum in the right-hand side is achieved at $\beta = a$, and

$$\lim_{n \to \infty} -\frac{1}{n} \log E[(\tilde{X}_n)^2] = a\alpha^* - H(\alpha^*) + L(a) < 2L(b).$$

Therefore, the traditional importance sampling scheme is not asymptotically optimal either.

TABLE 1
*Naive Monte Carlo scheme*

|                               | No. 1         | No. 2         | No. 3         | No. 4         |
|-------------------------------|---------------|---------------|---------------|---------------|
| Estimate $\hat{p}_n$ (%)      | 3.11          | 3.20          | 3.23          | 3.09          |
| Standard error (%)            | 0.17          | 0.18          | 0.18          | 0.17          |
| 95% confidence interval (%)   | [2.76, 3.46]  | [2.85, 3.55]  | [2.88, 3.58]  | [2.74, 3.44]  |

In Section 3 we argued the existence of asymptotically optimal adaptive importance sampling schemes in general. The construction of such adaptive schemes involved the selection of a nearly optimal control $\alpha^n = \{\alpha_j^n(\cdot,\cdot) : j = 0, 1, \ldots, n-1\}$. It was formally suggested in Section 4.1 that a good choice is to sample $\bar{Y}_j^n$, conditional on $\{\bar{Y}_i^n, i = 0, \ldots, j-1\}$, according to the transition probability matrix $Q_\alpha$ as in (4.4) with $\alpha$ being the conjugate of $\beta^*(x,t)$ given in (4.3), where $x = \bar{S}_j^n/n = (\bar{Y}_0^n + \cdots + \bar{Y}_{j-1}^n)/n$ and $t = 1 - j/n$. In case the conjugate of $\beta^*(x,t)$ is $\infty$ or $-\infty$, $\alpha$ is taken as a large positive or negative number; see Remark 4.3 for more details. The estimate is the sample mean of $K$ i.i.d. replications of

$$\bar{X}_n = \mathbb{1}_{\{\bar{S}_n^n/n \in A\}} \exp\left\{\sum_{j=1}^{n-1}(-\langle \alpha_j^n, g(\bar{Y}_j^n)\rangle + H(\alpha_j^n))\right\} \cdot \prod_{j=1}^{n-1} \frac{r(\bar{Y}_{j-1}^n; \alpha_j^n)}{r(\bar{Y}_j^n; \alpha_j^n)}.$$

The numerical results show that the controls constructed in this way have asymptotically optimal performance (Table 6).

The numerical results are reported in Tables 1–3 for $n = 60$. The theoretical value of $p_n$ is

$$p_n = P\{S_n/n \leq a\} + P\{S_n/n \geq b\}$$
$$= 0.83\% + 2.44\% = 3.27\%.$$

See Remark 4.2 for the computation of this theoretical value. For each table, we run four simulations each with sample size $K = 10{,}000$.

TABLE 2
*Traditional importance sampling scheme*

|                               | No. 1         | No. 2         | No. 3         | No. 4            |
|-------------------------------|---------------|---------------|---------------|------------------|
| Estimate $\hat{p}_n$ (%)      | 2.41          | 2.48          | 2.44          | 16.71            |
| Standard error (%)            | 0.04          | 0.04          | 0.04          | 14.22            |
| 95% confidence interval (%)   | [2.34, 2.48]  | [2.41, 2.56]  | [2.37, 2.51]  | [−11.73, 45.15]  |

TABLE 3
*Adaptive importance sampling scheme*

|  | No. 1 | No. 2 | No. 3 | No. 4 |
| --- | --- | --- | --- | --- |
| Estimate $\hat{p}_n$ (%) | 3.17 | 3.21 | 3.35 | 3.33 |
| Standard error (%) | 0.15 | 0.13 | 0.17 | 0.18 |
| 95% confidence interval (%) | [2.86, 3.47] | [2.85, 3.47] | [3.00, 3.69] | [2.96, 3.70] |

An interesting observation is that the traditional importance sampling scheme exhibits seemingly bizarre and inconsistent simulation results (Table 2). Similar phenomenon also occurs in the setting of Cramér's theorem, that is, where the Markov chain $Y$ reduces to a sequence of i.i.d. random variables; see [16, 19]. The explanation is also very similar. Under the alternative sampling distribution $Q_{\alpha^*}$, most of the sample means $\tilde{S}_n/n$ will end up near the point $b$. However, a few samples ("rogue" trajectories) have means that fall into the interval $(-\infty, a]$. Even though the "rogue" trajectories are rare, the Radon–Nikodym derivatives associated with them are so large that they dominate the variance. In simulation No. 4, the presence of a single "rogue" trajectory greatly raises the standard error associated with the estimate. Indeed, the proportion of the contribution to the second moment from this single "rogue" trajectories is more than 99%. In simulations No. 1, No. 2 and No. 3, however, there are no "rogue" trajectories, and the standard error associated with the estimate is deceptively small. The reason is that the standard error is itself estimated from the sample. Without "rogue" trajectories, we actually underestimate the standard error. Therefore, we cannot put much confidence in the standard errors thus obtained or in the "tight" confidence intervals that follow. Indeed, the confidence intervals from these three simulations do *not* contain the true value.

In contrast, the adaptive importance sampling, on the other hand, yields more accurate estimates and its performance is much more stable. Even though it does not show great advantage over naive Monte Carlo simulation for $n = 60$, it quickly does so when $n$ gets larger. The numerical results for different $n$ (with $K = 10{,}000$ fixed as before) are reported in Tables 4–6.

The naive Monte Carlo does not work well for bigger $n$. For $n = 120$ and $n = 180$, it yields estimates with large standard errors, and for $n = 240$, the simulation yields an estimate 0, that is, no sample mean reaches the target set $A$. As for the traditional importance sampling, each simulation gives a very "tight" confidence interval, due to the absence of "rogue" trajectories. However, as discussed before, we cannot put much belief into these estimates. Indeed, none of these confidence intervals cover the true value of $p_n$.

On the other hand, the adaptive importance scheme yields much more accurate estimates. In Table 6, the variable $\hat{V}^n$ denotes the sample estimate of the second moment $E[(\bar{X}_n)^2]$. Observe that as $n$ gets larger, the ratio

$$\frac{-(1/n)\log E[(\bar{X}_n)^2]}{-(1/n)\log p_n} = \frac{-\log E[(\bar{X}_n)^2]}{-\log p_n} \approx \frac{-\log \hat{V}^n}{-\log \hat{p}^n}$$

approaches 2. In other words, the adaptive importance sampling scheme is approaching optimality.

REMARK 4.2. The theoretical value of $p_n$ can be computed as follows. Let $X_n$ be the number of $-1$'s in a trajectory, that is,

$$X_n \doteq \sum_{j=0}^{n-1} \mathbb{1}_{\{Y_j=-1\}}.$$

Since $Y_0 \equiv 1$ and $Q(-1,1) = 1$, we have $0 \leq X_n \leq n/2$ with probability one. Clearly,

$$S_n = n - 2X_n,$$

whence it suffices to compute $P(X_n \geq m)$ for all nonnegative integers $m$ such that $2m \leq n$. But if we define $T_1 \doteq \inf\{j \geq 0 : Y_j = -1\}$, then $T_1 \geq 1$ and $T_1 - 1$ is geometrically distributed with parameter $(1-p)$. Moreover, $Y_{T_1} = -1$, and $Y_{1+T_1} = 1$. Now recursively define for $i \geq 2$, $T_i \doteq \inf\{j \geq 1 + T_{i-1} : Y_j = -1\}$. Then $\{T_1 - 1, T_2 - T_1 - 2, T_3 - T_2 - 2, \ldots\}$ is clearly a sequence of i.i.d. geometrically distributed random variables with parameter $(1-p)$, and

$$P(X_n \geq m) = P(T_m \leq n)$$
$$= P(T_m - 2m + 1 \leq n - 2m + 1).$$

But

$$T_m - 2m + 1 = (T_1 - 1) + (T_2 - T_1 - 2) + \cdots + (T_m - T_{m-1} - 2)$$

is the sum of i.i.d. geometrically distributed random variables, whence has a negative binomial distribution with parameter $m$ and $(1-p)$. Standard softwares such as SPLUS contain the cumulative distribution functions of negative binomial distributions, and can easily yield the desired probabilities.

TABLE 4
*Naive Monte Carlo simulation*

|  | $n=120$ | $n=180$ | $n=240$ |
|---|---|---|---|
| Theoretical $p_n$ | $1.61 \times 10^{-3}$ | $9.66 \times 10^{-5}$ | $6.35 \times 10^{-6}$ |
| Estimate $\hat{p}_n$ | $1.80 \times 10^{-3}$ | $20.00 \times 10^{-5}$ | 0 |
| Standard error | $0.42 \times 10^{-3}$ | $14.14 \times 10^{-5}$ | NA |
| 95% confidence interval | $[0.95, 2.65] \times 10^{-3}$ | $[-8.28, 48.28] \times 10^{-5}$ | NA |

TABLE 5
*Traditional importance sampling scheme*

|                          | $n = 120$                  | $n = 180$                   | $n = 240$                   |
|--------------------------|----------------------------|-----------------------------|-----------------------------|
| Theoretical $p_n$        | $1.61 \times 10^{-3}$      | $9.66 \times 10^{-5}$       | $6.35 \times 10^{-6}$       |
| Estimate $\hat{p}_n$     | $1.40 \times 10^{-3}$      | $8.76 \times 10^{-5}$       | $6.01 \times 10^{-6}$       |
| Standard error           | $0.02 \times 10^{-3}$      | $0.18 \times 10^{-5}$       | $0.13 \times 10^{-6}$       |
| 95% confidence interval  | $[1.35, 1.45] \times 10^{-3}$ | $[8.41, 9.12] \times 10^{-5}$ | $[5.74, 6.28] \times 10^{-6}$ |

REMARK 4.3. If $\beta^*(x,t) \geq 1$, then its conjugate is $\alpha^*(x,t) = +\infty$, in the sense that

$$L(\beta^*(x,t)) = \sup_\alpha [\alpha \beta^*(x,t) - H(\alpha)] = \lim_{\alpha \to +\infty} [\alpha \beta^*(x,t) - H(\alpha)].$$

The corresponding change of measure (at least formally) is

$$Q_{+\infty} \doteq \lim_{\alpha \to +\infty} Q_\alpha = \begin{bmatrix} 1 & 0 \\ 1 & 0 \end{bmatrix}.$$

Similarly, if $\beta^*(x,t) \leq 0$, then its conjugate is $\alpha^*(x,t) = -\infty$, with the corresponding change of measure

$$Q_{-\infty} \doteq \lim_{\alpha \to -\infty} Q_\alpha = \begin{bmatrix} 0 & 1 \\ 1 & 0 \end{bmatrix}.$$

However, neither of these two probability transition kernels is suitable for the purpose of importance sampling, since the probability measure induced by the original probability transition kernel $Q$ is *not* absolutely continuous with respect to the probability measure induced by $Q_{+\infty}$ or $Q_{-\infty}$.

To overcome this difficulty, we just take $\alpha$ to be a large positive or negative number whenever $\alpha^*(x,t) = +\infty$ or $\alpha^*(x,t) = -\infty$. In our numerical simulation, $\alpha$ is taken to be 5 if $\alpha^*(x,t) = +\infty$ and $-5$ if $\alpha^*(x,t) = -\infty$. The probability transition kernels corresponding to $\alpha = \pm 5$ are

$$Q_{+5} = \begin{bmatrix} 0.9999 & 0.0001 \\ 1 & 0 \end{bmatrix}, \qquad Q_{-5} = \begin{bmatrix} 0.0047 & 0.9953 \\ 1 & 0 \end{bmatrix},$$

TABLE 6
*Adaptive importance sampling scheme: asymptotic optimality*

|                          | $n = 120$                  | $n = 180$                   | $n = 240$                   |
|--------------------------|----------------------------|-----------------------------|-----------------------------|
| Theoretical $p_n$        | $1.61 \times 10^{-3}$      | $9.66 \times 10^{-5}$       | $6.35 \times 10^{-6}$       |
| Estimate $\hat{p}_n$     | $1.56 \times 10^{-3}$      | $9.73 \times 10^{-5}$       | $6.29 \times 10^{-6}$       |
| Standard error           | $0.04 \times 10^{-3}$      | $0.15 \times 10^{-5}$       | $0.07 \times 10^{-6}$       |
| 95% confidence interval  | $[1.49, 1.63] \times 10^{-3}$ | $[9.44, 10.02] \times 10^{-6}$ | $[6.15, 6.43] \times 10^{-6}$ |
| $(-\log \hat{V}^n)/(-\log \hat{p}_n)$ | 1.72          | 1.87                        | 1.93                        |

TABLE 7
*Traditional importance sampling scheme*

|  | No. 1 | No. 2 | No. 3 | No. 4 |
|---|---|---|---|---|
| Estimate $\hat{p}_n$ ($\times 10^{-5}$) | 2.14 | 2.37 | 2.29 | 9.20 |
| Standard error ($\times 10^{-5}$) | 0.11 | 0.15 | 0.14 | 6.85 |
| 95% confidence interval ($\times 10^{-5}$) | [1.92, 2.36] | [2.07, 2.67] | [2.01, 2.57] | [−4.50, 22.90] |

which are very close to $Q_{\pm\infty}$.

EXAMPLE 4.2. Consider a two-node tandem Jackson network with arrival rate $\lambda$ and consecutive service rates $\mu_1, \mu_2$. We assume the queueing system is stable, that is, $\lambda < \min\{\mu_1, \mu_2\}$, and, without loss of generality, $\lambda + \mu_1 + \mu_2 = 1$. The sizes of the first buffer and the second buffer are denoted by $B_1$ and $B_2$, respectively. Both buffer sizes are assumed to be finite.

We will work with the embedded discrete-time Markov chain $Y = \{Y_i = (Y_i^1, Y_i^2) : i = 0, 1, \ldots\}$, representing the queue lengths of the nodes at the epochs of transitions in the network. The chain $Y$ is irreducible and with finite state space $\mathcal{S} = \{(y_1, y_2) : y_i = 0, 1, \ldots, B_i; i = 1, 2\}$, whence uniformly recurrent. It is assumed throughout this example that the initial state is $Y_0 = (0, 0)$.

We are interested in estimating a class of probabilities associated with buffer overflow. More precisely, define $g = (g_1, g_2) : \mathcal{S} \to \{0, 1\}^2$ by

$$g_1(y) \doteq \mathbb{1}_{\{y_1 = B_1\}}, \qquad g_2(y) \doteq \mathbb{1}_{\{y_2 = B_2\}}$$

for every $y = (y_1, y_2) \in \mathcal{S}$, and let $S_n \doteq g(Y_0) + g(Y_1) + \cdots + g(Y_{n-1})$. We wish to estimate $p_n \doteq P\{S_n/n \in A\}$ for some Borel set $A$ of form

$$A = \{(x_1, x_2) : x_1 \geq \varepsilon_1 \text{ or } x_2 \geq \varepsilon_2\} \subset \mathbb{R}^2,$$

where $0 \leq \varepsilon_1, \varepsilon_2 \leq 1$. Note that the set $A$ is nonconvex.

The construction of the traditional and adaptive importance sampling schemes are very similar to Example 4.1. However, here the function $H : \mathbb{R}^2 \to \mathbb{R}$ and its conjugate $L : \mathbb{R}^2 \to \mathbb{R}^+$ do not admit closed-form expressions, and are computed numerically.

TABLE 8
*Adaptive importance sampling scheme*

|  | No. 1 | No. 2 | No. 3 | No. 4 |
|---|---|---|---|---|
| Estimate $\hat{p}_n$ ($\times 10^{-5}$) | 3.96 | 3.93 | 4.18 | 4.16 |
| Standard error ($\times 10^{-5}$) | 0.17 | 0.15 | 0.30 | 0.16 |
| 95% confidence interval ($\times 10^{-5}$) | [3.62, 4.30] | [3.63, 4.23] | [3.58, 4.78] | [3.84, 4.48] |

Analogously to Example 4.1, if we let $\beta^*$ be the minimizer that attains $\inf\{L(\beta): \beta \in A\}$ and let $\tilde{X}_n$ denote the traditional importance sampling estimate, then we have

$$(4.7) \qquad \lim_{n \to \infty} -\frac{1}{n} \log E[(\tilde{X}_n)^2] = \inf_{\beta \in A}[\alpha^*\beta - H(\alpha^*) + L(\beta)],$$

where $\alpha^*$ is the conjugate of $\beta^*$. It is not difficult to see that the traditional importance sampling scheme is asymptotically optimal if and only if $\beta^*$ is also a minimizer to the right-hand side of (4.7). However, this is often *not* the case, due to the nonconvexity of set $A$; see [16] for more discussion on this issue.

The simulation results for the traditional and adaptive schemes are reported in Tables 7 and 8. For comparison, the theoretical value of $p_n$ is also obtained via recursively computing the conditional distribution of $g(Y_k) + g(Y_{k+1}) + \cdots + g(Y_{n-1})$ given $Y_k$, for each $k = n-1, n-2, \ldots, 0$. Unlike Example 4.1, we choose not to report the results from naive Monte Carlo simulation (which is not asymptotically optimal). Actually, the naive Monte Carlo simulation, often giving an estimate 0 or an estimate with intolerably large standard error, is far inferior to either of the importance sampling schemes.

We choose $B_1 = B_2 = 6$, and $\lambda = 0.2$, $\mu_1 = \mu_2 = 0.4$. The state space $\mathcal{S}$ consists of $(B_1 + 1)(B_2 + 1) = 49$ states. Set $n = 50$ and $\varepsilon_1 = 0.3$, $\varepsilon_2 = 0.4$. Analogously to Example 4.1, one can check that the traditional importance sampling is not asymptotically optimal. Indeed, the infimum of $L(\beta)$ over set $A$ is attained at $\beta^* \approx (0.02, 0.4)$, while the minimizer for the right-hand side of (4.7) is $\bar{\beta} \approx (0.3, 0.01)$.

Each table consists of four simulation runs each with sample size $K = 10{,}000$. The theoretical value is $p_n = 4.10 \times 10^{-5}$.

The explanation for the behavior of traditional importance sampling (Table 8) is quite similar to that of Example 4.1—most of the sample means will end up near point $\beta^*$, while a few "rogue" trajectories will have means near point $\bar{\beta}$. Even though these "rogue" trajectories are rare, they carry huge Radon–Nikodym derivatives. Without the presence of "rogue" trajectories (simulations No. 1, No. 2 and No. 3), we have tight confidence intervals that we cannot put much faith in. With the presence of "rogue" trajectories (simulations No. 4), we get an estimate with very large standard error. On the contrast, the performance of adaptive schemes is much more stable and much better.

Similar phenomenon is also observed for various sets of parameters. We just list some numerical results in Tables 9 and 10 for the same setup, except the arrival rate and service rates are now $(\lambda, \mu_1, \mu_2) = (0.1, 0.4, 0.5)$. The sample size $K = 10{,}000$ is fixed as before. The erratic behavior of traditional schemes is more conspicuous. The asymptotic optimality of adaptive schemes is also clear from these numerical results.

TABLE 9
*Traditional importance sampling scheme*

|  | $n = 50$ | $n = 80$ | $n = 110$ |
|---|---|---|---|
| Theoretical $p_n$ | $5.15 \times 10^{-9}$ | $3.47 \times 10^{-12}$ | $1.83 \times 10^{-15}$ |
| Estimate $\hat{p}_n$ | $0.83 \times 10^{-9}$ | $0.81 \times 10^{-12}$ | $0.53 \times 10^{-15}$ |
| Standard error | $0.03 \times 10^{-9}$ | $0.02 \times 10^{-12}$ | $0.01 \times 10^{-15}$ |
| 95% confidence interval | $[0.77, 0.89] \times 10^{-9}$ | $[0.77, 0.85] \times 10^{-12}$ | $[0.51, 0.55] \times 10^{-15}$ |

## APPENDIX A

PROOF OF THEOREM 3.1. Proposition 2.1 and Condition 2.2 imply that

$$\lim_{n \to \infty} \frac{1}{n} \log P\{S_n/n \in A\} = -\inf_{\beta \in A} L(\beta).$$

Thanks to the discussion in Section 2.3, it suffices to show the lower bound (3.1), or

$$\liminf_n W^n \geq 2 \inf_{\beta \in A} L(\beta).$$

To this end, we extend the dynamics as in Section 3, and consider a mollified version of the original control problem. In other words, let $F: \mathbb{R}^d \to \mathbb{R}$ be an arbitrary bounded and Lipschitz continuous function, and define $V_F^n, W_F^n$ correspondingly; see the discussion from (3.1) to (3.2).

Since $V_F^n$ is the value function of a control problem, it satisfies the Bellman equation [4]

$$V_F^n(x, y; i)$$
$$= \inf_{\alpha \in \mathbb{R}^d} \int_{\mathcal{S}} e^{-2\langle \alpha, g(z) \rangle + 2H(\alpha)} \cdot \frac{r^2(y; \alpha)}{r^2(z; \alpha)} V_F^n\left(x + \frac{1}{n}g(z), z; i+1\right)$$
$$\times e^{\langle \alpha, g(z) \rangle - H(\alpha)} \cdot \frac{r(z; \alpha)}{r(y; \alpha)} p(y, dz)$$

TABLE 10
*Adaptive importance sampling scheme: asymptotic optimality*

|  | $n = 50$ | $n = 80$ | $n = 110$ |
|---|---|---|---|
| Theoretical $p_n$ | $5.15 \times 10^{-9}$ | $3.47 \times 10^{-12}$ | $1.83 \times 10^{-15}$ |
| Estimate $\hat{p}_n$ | $4.82 \times 10^{-9}$ | $3.36 \times 10^{-12}$ | $1.76 \times 10^{-15}$ |
| Standard error | $0.18 \times 10^{-9}$ | $0.11 \times 10^{-12}$ | $0.07 \times 10^{-15}$ |
| 95% confidence interval | $[4.46, 5.18] \times 10^{-9}$ | $[3.14, 3.58] \times 10^{-12}$ | $[1.62, 1.90] \times 10^{-6}$ |
| $(-\log \hat{V}^n)/(-\log \hat{p}_n)$ | 1.86 | 1.91 | 1.92 |

$$= \inf_{\alpha \in \mathbb{R}^d} \int_{\mathcal{S}} e^{-\langle \alpha, g(z) \rangle + H(\alpha)} \cdot \frac{r(y;\alpha)}{r(z;\alpha)} V_F^n\left(x + \frac{1}{n}g(y), z; i+1\right) p(y, dz),$$

together with terminal condition

$$V_F^n(x, y; n) = \exp\{-2nF(x)\}.$$

It follows from (3.2) that

(A.1)
$$W_F^n(x, y; i) = -\frac{1}{n} \log \inf_{\alpha \in \mathbb{R}^d} \int_{\mathcal{S}} e^{-\langle \alpha, g(z) \rangle + H(\alpha)}$$
$$\times \frac{r(y;\alpha)}{r(z;\alpha)} e^{-nW_F^n(x+1/ng(y), z; i+1)} p(y, dz)$$

and that $W_F^n(x, y; n) = 2F(x)$.

The discussion in Section 3 now prompts the following definition. Fixing an arbitrary $m \in \mathbb{N}$, for $0 \leq k \leq m - 1$, define recursively

(A.2)
$$U_F^m(x; k) = \sup_{\alpha \in \mathbb{R}^d} \inf_{\beta \in \mathbb{R}^d} \left[ U_F^m\left(x + \frac{1}{m}\beta; k+1\right) + \frac{1}{m}(L(\beta) + \langle \alpha, \beta \rangle - H(\alpha)) \right],$$

given the terminal condition

(A.3) $$U_F^m(x; m) \doteq 2F(x) \qquad \forall x \in \mathbb{R}^d.$$

See Section 3 for the interpretation of $W_F^n$ and $U_F^m$ as lower values of games. The key observation is the following lemma, whose proof is deferred to Appendix C.

LEMMA A.1. *For an arbitrary sequence $x^n \to x \in \mathbb{R}^d$, we have*

$$\liminf_{n \to \infty} \inf_{y \in \mathcal{S}} W_F^n(x^n, y; \lfloor nk/m \rfloor) \geq U_F^m(x; k), \qquad k = 0, 1, \ldots, m.$$

Assume Lemma A.1 holds for the moment. All that remains to show is the inequality

(A.4) $$\liminf_{m \to \infty} U_F^m(x; 0) \geq 2 \inf_{\beta \in \mathbb{R}^d} \{L(\beta) + F(x + \beta)\}.$$

Indeed, suppose (A.4) is true. Fix an arbitrary $j \in \mathbb{N}$, and define $F_j(y) \doteq j(d(y, \bar{A}) \wedge 1)$, which is bounded and Lipschitz continuous. Since $\mathbb{1}_A(y) \leq \exp\{-2nF_j(y)\}$, we have

$$\liminf_{n \to \infty} W^n \geq \liminf_{n \to \infty} W_{F_j}^n(0, y_0; 0)$$
$$\geq \liminf_{m \to \infty} U_{F_j}^m(0; 0)$$
$$\geq 2 \inf_{\beta \in \mathbb{R}^d} [L(\beta) + F_j(\beta)].$$

Exactly as in [14], pages 10 and 11, a compactness argument shows that

$$\lim_{j\to\infty} \inf_{\beta\in\mathbb{R}^d} \{L(\beta) + F_j(\beta)\} = \inf_{\beta\in\bar{A}} L(\beta),$$

and we complete the proof.

Now we show inequality (A.4). The idea is to represent $U_F^m$ as the value function of a control problem with the help of the min/max theorem. To this end, define

$$\mathcal{C} \doteq \left\{\theta \in \mathcal{P}(\mathbb{R}^d) : \int L(\beta)\theta(d\beta) < \infty\right\}$$

and rewrite (A.2) as

$$U_F^m(x;k) = \sup_{\alpha\in\mathbb{R}^d} \inf_{\theta\in\mathcal{C}} \left[\int U_F^m\left(x + \frac{1}{m}\beta; k+1\right)\theta(d\beta) \right.$$
$$\left. + \frac{1}{m}\left(\int L(\beta)\theta(d\beta) + \left\langle \alpha, \int \beta\theta(d\beta)\right\rangle - H(\alpha)\right)\right].$$

We make the following useful observation, whose proof is deferred to Appendix C.

LEMMA A.2. $U_F^m(\cdot;k)$ *is bounded and Lipschitz continuous for every $k$. Indeed,*

$$\|U_F^m(x;k)\| \leq 2\|F\|_\infty \qquad \forall x \in \mathbb{R}^d, \ k=0,1,\ldots,m,$$

*and $U_F^m(\cdot;k)$ is Lipschitz continuous with Lipschitz constant $2L_F$, where $L_F$ is the Lipschitz constant for the mollifier $F$.*

The next lemma is a version of min/max theorem, whose proof is almost identical to [16], Lemma 2.2, and thus omitted.

LEMMA A.3. *For any bounded and lower semicontinuous function $f:\mathbb{R}^d \to \mathbb{R}$, we have*

$$\sup_{\alpha\in\mathbb{R}^d} \inf_{\theta\in\mathcal{C}} \left[\int f(\beta)\,d\theta + \int L(\beta)\,d\theta + \left\langle\alpha, \int \beta\,d\theta\right\rangle - H(\alpha)\right]$$
$$= \inf_{\theta\in\mathcal{C}} \sup_{\alpha\in\mathbb{R}^d} \left[\int f(\beta)\,d\theta + \int L(\beta)\,d\theta + \left\langle\alpha, \int \beta\,d\theta\right\rangle - H(\alpha)\right].$$

Thanks to Lemmas A.2 and A.3, we obtain

$$U_F^m(x;k) = \inf_{\theta\in\mathcal{C}} \sup_{\alpha\in\mathbb{R}^d} \left[\int U_F^m\left(x + \frac{1}{m}\beta; k+1\right)\theta(d\beta)\right.$$

$$
\text{(A.5)} \quad + \frac{1}{m}\left(\int L(\beta)\theta(d\beta) + \left\langle \alpha, \int \beta\theta(d\beta)\right\rangle - H(\alpha)\right)\right]
$$

$$
= \inf_{\theta \in \mathcal{C}}\left[\int U_F^m\left(x + \frac{1}{m}\beta; k+1\right)\theta(d\beta)\right.
$$

$$
\left. + \frac{1}{m}\left(\int L(\beta)\theta(d\beta) + L\left(\int \beta\theta(d\beta)\right)\right)\right].
$$

This last display implies that $U_F^m$ has an interpretation as the minimal cost of a stochastic control problem. To simplify the notation, we state the control problem only for the case $k = 0$. The control problem will be defined on a probability $(\tilde{\Omega}, \tilde{\mathcal{F}}, \tilde{P})$, and $\tilde{E}_x$ will denote that the initial condition of the state process is $x$. An admissible control is a sequence $\{\nu_j^m, j = 0, 1, \ldots, m-1\}$, with each $\nu_j^m$ being a stochastic kernel on $\mathbb{R}^d$ given $\mathbb{R}^d$. Given an admissible control sequence, the state dynamics are defined by $\tilde{S}_0^m = mx$ and

$$
\tilde{S}_{j+1}^m \doteq \tilde{S}_j^m + \tilde{Y}_j^m,
$$

where

$$
\tilde{P}\{\tilde{Y}_j^m \in dy | \tilde{Y}_i^m, 0 \le i < j\} = \tilde{P}\{\tilde{Y}_j^m \in dy | \tilde{S}_j^m / m\} = \nu_j^m(dy | \tilde{S}_j^m / m).
$$

We then define the value function

$$
\tilde{v}_F^m(x; 0) \doteq \inf_{\{\nu_j^m\}} \tilde{E}_x\left[\sum_{j=0}^{m-1} \frac{1}{m}\left[\int L(y)\nu_j^m(dy) + L\left(\int y\nu_j^m(dy)\right)\right] + 2F(\tilde{S}_m^m/m)\right],
$$

where the infimum is taken over all controls $\{\nu_j^m\}$ and resulting controlled processes $\{\tilde{S}_j^m/m\}$ that start at $x$ at time 0. Since $\tilde{v}_F^m$ also satisfies the DPE (A.5) ([4], Chapter 8) and terminal condition $\tilde{v}_F^m(x; m) = U_F^m(x; m) = 2F(x)$, we obtain by induction that $U_F^m(x; k) = \tilde{v}_F^m(x; k)$ for all $x \in \mathbb{R}^d$ and $k \in \{0, \ldots, m\}$.

Define a stochastic kernel $\nu^m$ on $\mathbb{R}^d$ given $[0, 1]$ by

$$
\nu^m(dy|t) \doteq \begin{cases} \nu_j^m(dy), & \text{if } t \in [j/m, (j+1)/m), \ j = 0, 1, \ldots, m-2, \\ \nu_{m-1}^m(dy), & t \in [(m-1)/m, 1]. \end{cases}
$$

Let $\lambda$ denote Lebesgue measure. Then the definition of $\nu^m(dy|t)$ and the convexity of $L$ imply that

$$
U_F^m(x; 0) = \inf_{\{\nu_j^m\}} \tilde{E}_x\left[\int_0^1 \int_{\mathbb{R}^d} L(y)\nu^m(dy|t)\,dt \right.
$$

$$
\left. + \sum_{j=0}^{m-1} \frac{1}{m} L\left(\int_{\mathbb{R}^d} y\nu_j^m(dy)\right) + 2F(\tilde{S}_m^m/m)\right]
$$

DYNAMIC IMPORTANCE SAMPLING 27$$\geq \inf_{\{\nu_j^m\}} \tilde{E}_x \bigg[ \int_0^1 \int_{\mathbb{R}^d} L(y) \nu^m(dy|t) \, dt$$
$$+ L\bigg( \sum_{j=0}^{m-1} \frac{1}{m} \int_{\mathbb{R}^d} y \nu_j^m(dy) \bigg) + 2F(\tilde{S}_m^m/m) \bigg]$$
$$= \inf_{\{\nu_j^m\}} \tilde{E}_x \bigg[ \int_{\mathbb{R}^d \times [0,1]} L(y) \nu^m(dy \times dt)$$
$$+ L\bigg( \int_{\mathbb{R}^d \times [0,1]} y \nu^m(dy \times dt) \bigg) + 2F(\tilde{S}_m^m/m) \bigg],$$

where $\nu^m(dy \times dt) \doteq \nu^m(dy|t) \, dt$. A straightforward weak convergence approach will be adopted to derive the desired inequality (A.4). Since the proof is essentially the same as [14], Theorem 5.3.5, we only give a sketch.

For each $\varepsilon > 0$, there exist a sequence of controls $\{\nu^m, m \in \mathbb{N}\}$ such that, for every $m$, we have

$$U_F^m(x;0) + \varepsilon \geq \tilde{E}_x \bigg[ \int_{\mathbb{R}^d \times [0,1]} L(y) \, d\nu^m + L\bigg( \int_{\mathbb{R}^d \times [0,1]} y \, d\nu^m \bigg) + 2F(\tilde{S}_m^m/m) \bigg].$$

Furthermore, since $L$ is nonnegative and $F$ is bounded, we have

$$\sup_{m \in \mathbb{N}} \tilde{E}_x \int_{\mathbb{R}^d \times [0,1]} L(y) \nu^m(dy \times dt) < \infty.$$

However, since function $L$ is superlinear (Proposition 2.1), it is not difficult to check that $\{\nu^m\}$ is uniformly integrable in the sense that

$$\lim_{C \to \infty} \sup_{m \in \mathbb{N}} \tilde{E}_x \int_{\{y : \|y\| > C\} \times [0,1]} \|y\| \nu^m(dy \times dt) = 0.$$

It follows from that proof of [14], Proposition 5.3.2, that $\{\nu^m\}$ is indeed tight. Therefore, we can extract a weakly convergent sub-subsequence, still denoted by $\{\nu^m\}$, such that $\nu^m \Rightarrow \nu$ for some stochastic kernel $\nu$ whose second marginal is Lebesgue measure ([14], Lemma 5.3.4). We utilize the Skorokhod representation [6], which allows us to assume (when calculating the limits of the integrals) that the convergence is actually w.p.1. It follows from the uniform integrability of $\{\nu^m\}$ and the proof of [14], Proposition 5.3.5, that

$$\int_{\mathbb{R}^d \times [0,1]} y \nu^m(dy \times dt) \xrightarrow{P} \int_{\mathbb{R}^d \times [0,1]} y \nu(dy \times dt)$$

and

$$\tilde{S}_m^m/m \xrightarrow{P} Z \doteq x + \int_{\mathbb{R}^d \times [0,1]} y \nu(dy \times dt).$$

Furthermore, it follows from the lower-semicontinuity and nonnegativity of $L$ ([14], Lemma A.3.12) that, with probability one,

$$\liminf_m \int_{\mathbb{R}^d \times [0,1]} L(y)\nu^m(dy \times dt) \geq \int_{\mathbb{R}^d \times [0,1]} L(y)\nu(dy \times dt).$$

Thanks to convexity of $L$ and Jensen's inequality, we have

$$\int_{\mathbb{R}^d \times [0,1]} L(y)\nu(dy \times dt) \geq L\bigg(\int_{\mathbb{R}^d \times [0,1]} y\nu(dy \times dt)\bigg).$$

By Fatou's lemma and the lower-semicontinuity of $L$ [28], we have

$$\liminf_n U_F^m(x;0) + \varepsilon \geq \tilde{E}_x\bigg[2L\bigg(\int_{\mathbb{R}^d \times [0,1]} y\nu(dy \times dt)\bigg) + 2F(Z)\bigg].$$

It is now trivial that the right-hand side of the last inequality is bounded below by

$$2 \inf_{\beta \in \mathbb{R}^d}[L(\beta) + F(x+\beta)].$$

Since $\varepsilon > 0$ is arbitrary, (A.4) follows readily, which completes the proof. □

## APPENDIX B

**A large deviation upper bound.** In this section we study a uniform large deviation principle upper bound, which is essential for proving the key Lemma A.1. We present a proof based on the weak convergence approach [14]. Alternatively, one can adapt the methodology in [13].

The following two lemmas will be useful:

LEMMA B.1. *Suppose $\mathcal{S}$ is a Polish space and $\mathcal{P}(\mathcal{S})$ is the space of probability measures on $\mathcal{S}$ endowed with the weak convergence topology. Consider a sequence of random variables $\mu_n : (\Omega^n, \mathcal{F}^n, P^n) \to \mathcal{P}(\mathcal{S})$. In other words, $\{\mu^n\}$ is a sequence of random probability measures. Then $\{\mu^n\}$ is tight if and only if the sequence $\{E^n \mu^n\}$ is tight. Here $E^n \mu^n \in \mathcal{P}(\mathcal{S})$ is defined by*

$$(E^n \mu^n)(A) \doteq \int_{\Omega^n} \mu^n(\omega)(A) P^n(d\omega)$$

*for every Borel set $A$ in $\mathcal{P}(\mathcal{S})$.*

PROOF. See [23], Theorem 6.1, Chapter 1. □

LEMMA B.2. *Suppose $\mathcal{S}$ is a Polish space, $\{\mu^n\} \subset \mathcal{P}(\mathcal{S})$, and $p(\cdot, \cdot)$ a probability transition kernel. If $\mu^n \to \mu$ in the $\tau$-topology for some $\mu \in \mathcal{P}(\mathcal{S})$, then*

$$\mu^n \otimes p \to \mu \otimes p$$

in the $\tau$-topology. Here $\mu \otimes p$ denotes the probability measure on $\mathcal{S} \times \mathcal{S}$ given by

$$(\mu \otimes p)(B) \doteq \int_B \mu(dx)p(x,dy)$$

for every Borel set $B \subseteq \mathcal{S} \times \mathcal{S}$.

PROOF. It suffices to show that

$$\int_{\mathcal{S} \times \mathcal{S}} f(x,y)\mu^n(dx)p(x,dy) \to \int_{\mathcal{S} \times \mathcal{S}} f(x,y)\mu(dx)p(x,dy)$$

for every bounded, measurable function $f$. Since $\mu^n \to \mu$ in the $\tau$-topology, it remains to show that

$$\int_{\mathcal{S}} f(x,y)p(x,dy)$$

is a bounded and measurable function (over $x$). The boundedness is trivial, and the measurability follows from Fubini's theorem; compare [5], Exercise 18.20. □

PROPOSITION B.3. *Suppose $Y = \{Y_j, j \in \mathbb{N}\}$ is a Markov chain that takes values in a Polish space $\mathcal{S}$. Let $p$ denote the probability transition kernel of $Y$, and assume Condition* 2.1 *holds. Suppose $g : \mathcal{S} \to \mathbb{R}^d$ is a bounded measurable function, and define $H$ and $L$ as in* (2.1)–(2.3). *Then for any fixed $\alpha \in \mathbb{R}^d$, bounded and continuous function $f : \mathbb{R}^d \to \mathbb{R}$, and sequence $x^n \to x \in \mathbb{R}^d$, we have*

$$\liminf_{n\to\infty} \inf_{y \in \mathcal{S}} -\frac{1}{n} \log E_y \Bigg[ \exp\bigg\{-\bigg\langle \alpha, \sum_{j=0}^{n-1} g(Y_j) \bigg\rangle \bigg\}$$
$$\times \exp\bigg\{-nf\bigg(x^n + \frac{1}{n}\sum_{j=0}^{n-1} g(Y_j)\bigg)\bigg\}\Bigg] \geq I_f(x),$$

*where*

$$I_f(x) \doteq \inf_\beta [f(x+\beta) + L(\beta) + \langle \alpha, \beta \rangle].$$

PROOF. Let

$$v^n(x,y) \doteq -\frac{1}{n} \log E_y \Bigg[ \exp\bigg\{-\bigg\langle \alpha, \sum_{j=0}^{n-1} g(Y_j) \bigg\rangle\bigg\} \exp\bigg\{-nf\bigg(x + \frac{1}{n}\sum_{j=0}^{n-1} g(Y_j)\bigg)\bigg\}\Bigg]$$

$$= -\frac{1}{n} \log \int \exp\bigg\{-\bigg\langle \alpha, \sum_{j=0}^{n-1} g(y_j) \bigg\rangle\bigg\} \exp\bigg\{-nf\bigg(x + \frac{1}{n}\sum_{j=0}^{n-1} g(y_j)\bigg)\bigg\} d\pi_y^n.$$

Here $\pi_y^n$ is the joint distribution of $(Y_0, Y_1, \ldots, Y_{n-1})$, or

$$\pi_y^n(dy_0, dy_1, \ldots, dy_n) \doteq \delta_y(dy_0) p(y_0, dy_1) p(y_1, dy_2) \cdots p(y_{n-1}, dy_n).$$

Clearly $v^n$ is bounded, thanks to the boundedness of $g$ and $f$. It suffices to show that for every sequence $x^n \to x$ and $\{y^n\} \subseteq \mathcal{S}$,

(B.1) $$\liminf_{n \to \infty} v^n(x^n, y^n) \geq I_f(x).$$

For an arbitrary $\varepsilon > 0$, the relative entropy representation of exponential integrals (3.3) ([14], Proposition 1.4.2) yields the existence of a probability measure $\mu^n$ on $\mathcal{S}^{n+1}$ such that

(B.2)
$$v^n(x^n, y^n) + \varepsilon \geq \frac{1}{n} R(\mu^n \| \pi_{y^n}^n) + \left\langle \alpha, \int \frac{1}{n} \sum_{j=0}^{n-1} g(y_j) \, d\mu^n \right\rangle$$
$$+ \int f\left(x^n + \frac{1}{n} \sum_{j=0}^{n-1} g(y_j)\right) d\mu^n.$$

In particular, it is not hard to see that

(B.3) $$\sup_{n \in \mathbb{N}} \frac{1}{n} R(\mu^n \| \pi_{y^n}^n) < \infty.$$

We can factor $\mu^n$ as in [14], Theorems A.5.4 and A.5.6:

$$\mu^n(dy_0, dy_1, \ldots, dy_n) = \mu_0^n(dy_0) \mu_1^n(dy_1|y_0) \cdots \mu_n^n(dy_n|y_{n-1}, y_{n-2}, \ldots, y_0).$$

Now consider a probability space $(\tilde{\Omega}, \tilde{\mathcal{F}}, \tilde{P})$, on which we define a stochastic process given by

$$\tilde{P}(\tilde{Y}_0^n \in dy) = \mu_0^n(dy_0),$$
$$\tilde{P}(\tilde{Y}_{j+1}^n \in dy | \tilde{Y}_i^n, i = 0, 1, \ldots, j) = \mu_{j+1}^n(dy | \tilde{Y}_i^n, i = 0, 1, \ldots, j)$$

for $j = 0, 1, \ldots, n-1$. To ease exposition, let

$$\bar{\mu}_{j+1}^n(dy) \doteq \mu_{j+1}^n(dy | \tilde{Y}_i^n, i = 0, 1, \ldots, j),$$

which is a random probability measure on $\mathcal{S}$. Also define a random probability measure on $\mathcal{S} \times \mathcal{S}$ by

$$\gamma^n(dx \times dy) = \frac{1}{n} \sum_{j=0}^{n-1} \delta_{\tilde{Y}_j^n}(dx) \times \bar{\mu}_{j+1}^n(dy),$$

whose marginals are

$$(\gamma^n)_1 = \frac{1}{n} \sum_{j=0}^{n-1} \delta_{\tilde{Y}_j^n} \doteq \tilde{L}^n, \qquad (\gamma^n)_2 = \frac{1}{n} \sum_{j=0}^{n-1} \bar{\mu}_{j+1}^n.$$

Thanks to the chain rule [14], Theorem B.2.1, we have

$$\text{(B.4)} \quad \frac{1}{n}R(\mu^n\|\pi^n_{y^n}) = \frac{1}{n}\tilde{E}\left[R(\mu^n_0\|\delta_{y^n}) + \sum_{j=0}^{n-1} R(\bar{\mu}^n_{j+1}(\cdot)\|p(\tilde{Y}_j,\cdot))\right].$$

However,

$$\frac{1}{n}\sum_{j=0}^{n-1} R(\bar{\mu}^n_{j+1}(\cdot)\|p(\tilde{Y}^n_j,\cdot))$$

$$= \frac{1}{n}\sum_{j=0}^{n-1} R(\delta_{\tilde{Y}^n_j}(dy) \times \bar{\mu}^n_{j+1}(dz) \| \delta_{\tilde{Y}^n_j}(dy) \times p(\tilde{Y}^n_j, dz))$$

$$= \frac{1}{n}\sum_{j=0}^{n-1} R(\delta_{\tilde{Y}^n_j}(dy) \times \bar{\mu}^n_{j+1}(dz) \| \delta_{\tilde{Y}^n_j}(dy) \otimes p(y, dz))$$

$$\geq R\left(\frac{1}{n}\sum_{j=0}^{n-1} \delta_{\tilde{Y}^n_j}(dy) \times \bar{\mu}^n_{j+1}(dz) \Big\| \frac{1}{n}\sum_{j=0}^{n-1} \delta_{\tilde{Y}^n_j}(dy) \otimes p(y, dz)\right)$$

$$= R(\gamma^n \| \tilde{L}^n \otimes p),$$

where the inequality follows from the convexity of relative entropy $R(\cdot\|\cdot)$. Thanks to (B.2), and observing that

$$\int \frac{1}{n}\sum_{j=0}^{n-1} g(y_j) \, d\mu^n = \tilde{E}\int g \, d\tilde{L}^n,$$

$$\int f\left(x^n + \frac{1}{n}\sum_{j=0}^{n-1} g(y_j)\right) d\mu^n = \tilde{E}f\left(x^n + \int g \, d\tilde{L}^n\right),$$

we arrive at

$$v^n(x^n, y^n) + \varepsilon \geq \tilde{E}\left[R(\gamma^n\|\tilde{L}^n \otimes p) + \left\langle \alpha, \int g \, d\tilde{L}^n \right\rangle + f\left(x^n + \int g \, d\tilde{L}^n\right)\right].$$

It suffices to show $\{\gamma^n\}$ is tight. Indeed, if this is true, the same argument as in [14], Theorem 8.2.8, allows us to extract a weak convergent subsequence of $(\gamma^n, \tilde{L}^n)$, still indexed by $n$, such that

$$(\gamma^n, \tilde{L}^n) \Rightarrow (\gamma, \tilde{L})$$

for some stochastic kernel $\gamma$ on $\mathcal{S} \times \mathcal{S}$ and some stochastic kernel $\tilde{L}$ on $\mathcal{S}$, and a (random) transition probability function $q$ such that

$$\gamma(dy \times dz) = \tilde{L}(dy) \otimes q(y, dz)$$

and

(B.5) $$\tilde{L}q = \tilde{L}$$

hold almost surely. In particular, we have

$$(\gamma^n)_2 \Rightarrow (\gamma)_2 = \tilde{L}q = \tilde{L}.$$

Note (B.5) says that $\tilde{L}$ is indeed the invariant measure for the transition probability function $q$. Also observe that

$$\sup_{n \in \mathbb{N}} \tilde{E} R(\gamma^n \| \tilde{L}^n \otimes p) < \infty.$$

This implies the existence of a subsequence, still indexed by $n$, such that

$$\tilde{L}^n \to \tilde{L}, \qquad (\gamma^n)_2 \to \tilde{L}$$

in the $\tau$-topology; see the proof of [14], Lemma 9.3.3. Therefore,

$$\int g \, d\tilde{L}^n \to \int g \, d\tilde{L}$$

almost surely. Furthermore, thanks to Lemma B.2, $\tilde{L}^n \otimes p \to \tilde{L} \otimes p$ in the $\tau$-topology (hence, in the weak-topology) almost surely. The lower semi-continuity of $R(\cdot\|\cdot)$ implies

$$\liminf_{n \to \infty} R(\gamma^n \| \tilde{L}^n \otimes p) \geq R(\gamma \| \tilde{L} \otimes p).$$

It follows readily from Fatou's lemma that

$$\liminf_{n \to \infty} v^n(x^n, y^n) + \varepsilon$$
$$\geq \tilde{E}\left[R(\gamma \| \tilde{L} \otimes p) + \left\langle \alpha, \int g \, d\tilde{L} \right\rangle + f\left(x + \int g \, d\tilde{L}\right)\right]$$
$$= \tilde{E}\left[R(\tilde{L} \otimes q \| \tilde{L} \otimes p) + \left\langle \alpha, \int g \, d\tilde{L} \right\rangle + f\left(x + \int g \, d\tilde{L}\right)\right]$$
$$\geq \inf_{\{\mu q = \mu\}} \left[R(\mu \otimes q \| \mu \otimes p) + \left\langle \alpha, \int g \, d\mu \right\rangle + f\left(x + \int g \, d\mu\right)\right].$$

Recalling (2.4) and letting $\varepsilon \to 0$, we obtain

$$\liminf_{n \to \infty} v^n(x^n, y^n) \geq \inf_\beta [L(\beta) + \langle \alpha, \beta \rangle + f(x + \beta)],$$

which is the desired inequality (B.1).

It remains to show the tightness of $\{\gamma^n\}$. All we need is the tightness of the two marginals, $\{(\gamma^n)_1\}$ and $\{(\gamma^n)_2\}$. However, it is not difficult to observe that

$$\tilde{E}(\gamma^n)_1 = \tilde{E}\tilde{L}^n = \frac{1}{n}\sum_{j=0}^{n-1} \tilde{E}\delta_{\tilde{Y}_j^n} = \frac{1}{n}\sum_{j=0}^{n-1} \mu^{n,j},$$

where $\mu^{n,j}$ denotes the $j$th marginal of the probability $\mu^n$ and, similarly,

$$\tilde{E}(\gamma^n)_2 = \frac{1}{n}\sum_{j=0}^{n-1} \tilde{E}\bar{\mu}_{j+1}^n = \frac{1}{n}\sum_{j=0}^{n-1} \mu^{n,j+1}.$$

Letting $\|\cdot\|_{\mathrm{v}}$ denote the total variation metric, we have

(B.6) $$\|\tilde{E}(\gamma^n)_1 - \tilde{E}(\gamma^n)_2\|_{\mathrm{v}} = \frac{1}{n}\|\mu^{n,0} - \mu^{n,n}\|_{\mathrm{v}} \leq \frac{2}{n}.$$

If we can show $\{(\gamma^n)_2\}$ is tight, then Lemma B.1 implies $\{\tilde{E}(\gamma^n)_2\}$ is tight, which in turns yields the tightness of $\{\tilde{E}(\gamma^n)_1\}$, thanks to (B.6). Applying Lemma B.1 once again, we have the tightness of $\{(\gamma^n)_1\}$. Therefore, it is sufficient to show that $\{(\gamma^n)_2\}$ is tight. The proof will distinguish two cases: $m_0 = 1$ and $m_0 > 1$.

Suppose that $m_0 = 1$. Note that the nonnegativity of relative entropy, (B.3) and (B.4) imply

$$\sup_{n\in\mathbb{N}} \tilde{E}\frac{1}{n}\sum_{j=0}^{n-1} R(\bar{\mu}_{j+1}^n(\cdot)\|p(\tilde{Y}_j^n,\cdot)) < \infty.$$

It follows from the assumption of uniform recurrency

$$a\nu_p(\cdot) \leq p(y,\cdot) \leq b\nu_p(\cdot) \qquad \forall\, y \in \mathcal{S},$$

that

$$R(\bar{\mu}_{j+1}^n(\cdot)\|p(\tilde{Y}_j^n,\cdot)) \geq cR(\bar{\mu}_{j+1}^n\|\nu_p)$$

for some constant $c > 0$. It is now easy to derive from the convexity of relative entropy that

$$\sup_{n\in\mathbb{N}} \tilde{E} R((\gamma^n)_2\|\nu_p) \leq \sup_{n\in\mathbb{N}} \tilde{E}\frac{1}{n}\sum_{j=0}^{n-1} R(\bar{\mu}_{j+1}^n\|\nu_p) < \infty,$$

which further implies the tightness of $\{(\gamma^n)_2\}$ since $R(\cdot\|\nu_p)$ is a tightness function on $\mathcal{P}(\mathcal{S})$. Note that

$$\tilde{E}(\gamma^n)_2 = \frac{1}{n}\sum_{j=0}^{n-1} \mu^{n,j+1}$$

is also tight, thanks to Lemma B.1.

The general case with $m_0 > 1$ is slightly more complicated. We will give a proof with $m_0 = 2$, and observe that the proof for $m_0 > 2$ is essentially the same and thus omitted. Without loss of generality, we show $\{(\gamma^n)_2 : n \text{ even}\}$ to be tight. The tightness for $\{(\gamma^n)_2 : n \text{ odd}\}$ is similar.

To ease notation, let $\pi^n \doteq \pi^n_{y^n}$, and $\pi^{n,e}$ be the marginal distribution of $\pi^n$ over even coordinates; that is,

$$\pi^{n,e}(dy_0, dy_2, \ldots, dy_{n-2}, dy_n) = \delta_y(dy_0) p^{(2)}(y_0, dy_2) \cdots p^{(2)}(y_{n-2}, dy_n).$$

One can similarly define $\mu^{n,e}$, or

$$\mu^{n,e}(dy_0, dy_2, \ldots, dy_n) = \mu_0^n(dy_0) \mu_2^{n,e}(dy_2|y_0) \cdots \mu_n^{n,e}(dy_n|y_{n-2}, \ldots, y_2, y_0).$$

Thanks to the chain rule ([14], Theorem B.2.1) and nonnegativity of the relative entropy, we have

$$R(\mu^{n,e} \| \pi^{n,e}) \leq R(\mu^n \| \pi^n),$$

and, thus, $\sup_n \frac{1}{n} R(\mu^{n,e} \| \pi^{n,e}) < \infty$. With the same proof as for the case $m_0 = 1$, we have that

$$\frac{2}{n} \sum_{j=0}^{(n/2)-1} \mu^{n,e,j+1}$$

is tight; here $\mu^{n,e,j}$ is the $j$th marginal of $\mu^{n,e}$; that is,

$$\mu^{n,e,j}(dy_{2j}) = \mu^{n,e}(\mathcal{S}, \ldots, \mathcal{S}, dy_{2j}, \mathcal{S}, \ldots, \mathcal{S}).$$

One can similarly define $\mu^{n,o}$ as the marginal distribution of $\mu^n$ over odd coordinates, and the same argument can be carried over to prove the tightness of

$$\frac{2}{n} \sum_{j=0}^{(n/2)-1} \mu^{n,o,j+1}.$$

However, observe that

$$\mu^{n,e,j} = \mu^{n,2j}, \qquad \mu^{n,o,j} = \mu^{n,2j+1}.$$

We have

$$\frac{1}{2} \left( \frac{2}{n} \sum_{j=0}^{(n/2)-1} \mu^{n,e,j+1} + \frac{2}{n} \sum_{j=0}^{(n/2)-1} \mu^{n,o,j+1} \right) = \frac{1}{2} \sum_{j=0}^{n-1} \mu^{n,j+1} = \tilde{E}(\gamma^n)_2.$$

This implies the tightness of $\{\tilde{E}(\gamma^n)_2 : n \text{ even}\}$, which is equivalent to the tightness of $\{(\gamma^n)_2 : n \text{ even}\}$, thanks to Lemma B.1. $\square$

# APPENDIX C

**Proofs of Lemmas A.1 and A.2.**

PROOF OF LEMMA A.2. That $U_F^m(\cdot;k)$ is Lipschitz continuous with Lipschitz constant $2L_F$ follows trivially by induction and the terminal condition (A.3).

As for the boundedness of $U_F^m(\cdot;k)$, we first show it is bounded from below. Since $H(0) = 0$ and $L$ is nonnegative, definition (A.2) gives

$$U_F^m(x;k) \geq \inf_{\beta \in \mathbb{R}^d} \left[ U_F^m\left(x + \frac{1}{m}\beta; k+1\right) + \frac{1}{m}L(\beta) \right]$$

$$\geq \inf_{\beta \in \mathbb{R}^d} U_F^m\left(x + \frac{1}{m}\beta; k+1\right)$$

$$= \inf_{z \in \mathbb{R}^d} U_F^m(z; k+1)$$

for every $x$. It follows that, for every $k$,

$$\inf_{x \in \mathbb{R}^d} U_F^m(x;k) \geq \inf_{x \in \mathbb{R}^d} U_F^m(x;k+1) \geq \cdots \geq \inf_{x \in \mathbb{R}^d} U_F^m(x;m) \geq -2\|F\|_\infty.$$

It remains to show that $U_F^m$ is bounded from above. Let $\bar{\beta}$ be a subdifferential of the convex function $H$ at $\alpha = 0$. Then

$$L(\bar{\beta}) = \sup_{\alpha \in \mathbb{R}^d} [\langle \alpha, \bar{\beta}\rangle - H(\alpha)] = 0$$

and the supremum is achieved at $\alpha = 0$. By definition (A.2) again, we have

$$U_F^m(x;k) \leq \sup_{\alpha \in \mathbb{R}^d} \left[ U_F^m\left(x + \frac{1}{m}\bar{\beta}; k+1\right) + \frac{1}{m}(\langle \alpha, \bar{\beta}\rangle - H(\alpha)) \right]$$

$$= U_F^m\left(x + \frac{1}{m}\bar{\beta}; k+1\right)$$

$$\leq \sup_{z \in \mathbb{R}^d} U_F^m(z; k+1)$$

for every $x$. It follows that, for every $k$,

$$\sup_{x \in \mathbb{R}^d} U_F^m(x;k) \leq \sup_{x \in \mathbb{R}^d} U_F^m(x;k+1) \leq \cdots \leq \sup_{x \in \mathbb{R}^d} U_F^m(x;m) \leq 2\|F\|_\infty.$$

This completes the proof. □

PROOF OF LEMMA A.1. The proof is by induction. For $k = m$, we have

$\lfloor nk/m \rfloor = n$. By definition,

$$\liminf_{n \to \infty} \inf_{y \in \mathcal{S}} W_F^n(x^n, y; n) = \liminf_{n \to \infty} \inf_{y \in \mathcal{S}} 2F(x^n) = \liminf_{n \to \infty} 2F(x^n),$$

and Lemma A.1 follows trivially from the continuity of $F$.

Assume now the claim holds for $k+1$. Let $\ell(n) \doteq \lfloor n(k+1)/m \rfloor - \lfloor nk/m \rfloor$. Also, let $\pi_y^j$ be the probability measure on $\mathcal{S}^{j+1}$ defined by

$$\pi_y^j(dy_0, dy_1, \ldots, dy_j) \doteq \delta_y(dy_0) p(y_0, dy_1) p(y_1, dy_2) \cdots p(y_{j-1}, dy_j)$$

for every $y \in \mathcal{S}$ and every $j \in \mathbb{N}$.

For an arbritrary $\alpha \in \mathbb{R}^d$, let

$$U_{\alpha,F}^m(x; k) \doteq \inf_{\beta \in \mathbb{R}^d} \left[ U_F^m\left(x + \frac{1}{m}\beta; k+1\right) + \frac{1}{m}(L(\beta) + \langle \alpha, \beta \rangle - H(\alpha)) \right].$$

It follows from the definition that $U_F^m(x; k) = \sup_\alpha U_{\alpha,F}^m(x; k)$. Therefore, all we need to show is that, for every $\alpha \in \mathbb{R}^d$ and any sequence $x^n \to x$,

$$\liminf_{n \to \infty} \inf_{y \in \mathcal{S}} W_F^n(x^n, y; \lfloor nk/m \rfloor) \geq U_{\alpha,F}^m(x; k).$$

However, for an arbitrary fixed $\alpha \in \mathbb{R}^d$, the dynamic programming principle implies that

$$W_F^n(x, y; \lfloor nk/m \rfloor)$$

$$\geq -\frac{1}{n} \log \int \exp\left\{ -\left\langle \alpha, \sum_{j=1}^{\ell(n)} g(y_i) \right\rangle + \ell(n) H(\alpha) \right\} \cdot \prod_{j=1}^{\ell(n)} \frac{r(y_{j-1}; \alpha)}{r(y_j; \alpha)}$$

$$\times \exp\left\{ -n W_F^n\left( x + \frac{1}{n} \sum_{j=0}^{\ell(n)-1} g(y_j), y_{\ell(n)}; \lfloor n(k+1)/m \rfloor \right) \right\} d\pi_y^{\ell(n)}$$

$$= -\frac{1}{n} \log \int \exp\left\{ -\left\langle \alpha, \sum_{j=1}^{\ell(n)} g(y_i) \right\rangle + \ell(n) H(\alpha) \right\} \cdot \frac{r(y_0; \alpha)}{r(y_{\ell(n)}; \alpha)}$$

$$\times \exp\left\{ -n W_F^n\left( x + \frac{1}{n} \sum_{j=0}^{\ell(n)-1} g(y_j), y_{\ell(n)}; \lfloor n(k+1)/m \rfloor \right) \right\} d\pi_y^{\ell(n)}.$$

Since $g$ is bounded and $r(\cdot;\alpha)$ is both bounded from above and bounded away from zero by (2.2), it suffices to show

(C.1) $$\liminf_{n\to\infty} \inf_{y\in\mathcal{S}} \bar{v}_F^n(x^n, y; 0) \geq U_{\alpha,F}^m(x; k),$$

where

$$\bar{v}_F^n(x, y; 0)$$
$$\doteq -\frac{1}{n}\log \int \exp\left\{-\left\langle \alpha, \sum_{j=0}^{\ell(n)-1} g(y_i) \right\rangle + \ell(n)H(\alpha)\right\}$$
$$\times \exp\left\{-nW_F^n\left(x + \frac{1}{n}\sum_{j=0}^{\ell(n)-1} g(y_j), y_{\ell(n)}; \lfloor n(k+1)/m \rfloor\right)\right\} d\pi_y^{\ell(n)}.$$

We claim that inequality (C.1) is a direct consequence of

(C.2) $$\liminf_{n\to\infty} \inf_{y\in\mathcal{S}} v_F^n(x^n, y; 0) \geq U_{\alpha,F}^m(x; k),$$

where

$$v_F^n(x, y; 0) \doteq -\frac{1}{n}\log \int \exp\left\{-\left\langle \alpha, \sum_{j=0}^{\ell(n)-1} g(y_i) \right\rangle + \ell(n)H(\alpha)\right\}$$
$$\times \exp\left\{-nU_F^m\left(x + \frac{1}{n}\sum_{j=0}^{\ell(n)-1} g(y_j); k+1\right)\right\} d\pi_y^{\ell(n)}.$$

Indeed, since $\ell(n) \leq n$, one can always find a compact set $K \subseteq \mathbb{R}^d$ such that

$$x^n + \frac{1}{n}\sum_{j=0}^{\ell(n)-1} g(y_j) \in K \qquad \forall (y_0, y_1, \ldots, y_{\ell(n)}), \ \forall n \in \mathbb{N},$$

thanks to the boundedness of $g$ and the assumption $x^n \to x$. It is also not hard to show by contradiction from the induction hypothesis and the continuity of $U_F^m$ (Lemma A.2) that, for any $\varepsilon > 0$, there exists $N(\varepsilon) \in \mathbb{N}$ such that for all $x \in K$ and $n \geq N(\varepsilon)$,

$$\inf_{y\in\mathcal{S}} W_F^n(x, y; \lfloor n(k+1)/m \rfloor) - U_F^m(x; k+1) \geq -\varepsilon.$$

We arrive at

$$\liminf_{n\to\infty} \inf_{y\in\mathcal{S}} \bar{v}_F^n(x^n, y; 0) \geq \liminf_{n\to\infty} \inf_{y\in\mathcal{S}} v_F^n(x^n, y; 0) - \varepsilon$$

for every $\varepsilon > 0$. It follows that (C.1) is implied by (C.2).

It remains to show (C.2), which is an easy consequence of the uniform large deviation bound Proposition B.3, Lemma A.2, boundedness of $g$, and that

$$\left|\frac{\ell(n)}{n} - \frac{1}{m}\right| \to 0.$$

This completes the proof. □

**Acknowledgments.** We are indebted to Associate Editor and the referee for their careful reading of the first version and many helpful suggestions. In particular, Example 4.2 is inspired by the comments from Associate Editor.

DIVISION OF APPLIED MATHEMATICS
BROWN UNIVERSITY
PROVIDENCE, RHODE ISLAND 02912
USA
E-MAIL: dupuis@dam.brown.edu
E-MAIL: huiwang@cfm.brown.edu